\documentclass{article}

\usepackage{amssymb}
\usepackage{latexsym}
\usepackage[mathscr]{euscript}
\usepackage{pb-diagram}
\usepackage{latexsym}
\usepackage{epsfig}

\def\column#1#2{\mathrel{\mathop{#1}\limits_{#2}}}

\newcommand\ul {\underline}

\newcommand\commentout[1]{\marginpar{\tiny $\backslash$commentout}}

\newcommand\qed{\hfill$\square$}
\newcommand \kay{K}

\newcommand \comp[2]{ {#1}^{\wedge}_{#2}}

\newcommand \homo{homomorphism}

\newtheorem{Lemma}{Lemma}[section]
 
\newtheorem{Proposition}[Lemma]{Proposition}

\newtheorem{Remark}[Lemma]{Remark}

\newtheorem{gLemma}{Lemma}[section]
\newtheorem{gTheorem}[gLemma]{Theorem}
\newtheorem{gProposition}[gLemma]{Proposition}
\newtheorem{gDefinition}[gLemma]{Definition}
\newtheorem{gCorollary}[gLemma]{Corollary}
\newtheorem{gConjecture}[gLemma]{Conjecture}

\newtheorem{GLemma}{Lemma}[section]
\newtheorem{GTheorem}[GLemma]{Theorem}
\newtheorem{GProposition}[GLemma]{Proposition}

\newenvironment{Proof}{\par\noindent\textbf{Proof:}}
{\qed}

\usepackage{graphics}

\begin{document}

\begin{center}
{\bf  \Large Derived Representation Theory and 
the Algebraic $K$-theory of Fields}
\end{center}

\begin{center}
{\bf Gunnar Carlsson\footnote{Supported in part by NSF DMS 0104162} \\Department of Mathematics
\\Stanford University 
\\ Stanford, California 94305}
\end{center}

\section{Introduction}
\label{Introduction}
 Quillen's higher algebraic $K$-theory for fields $F$ has been the object of intense study since its introduction in 1972 \cite{Quillen}.  The main direction of research has been the construction of  ``descent spectral sequences" whose $E_2$-term involved the cohomology of the absolute Galois group $G_F$ with coefficients in the algebraic $K$--theory of an algebraically closed field.  The form of such a spectral sequence was conjectured in \cite{Quillen3} and \cite{Lichtenbaum}, but it was soon realized that it could not be expected to converge  to algebraic $K$-theory exactly.  However,  it appeared likely that it could converge to algebraic $K$-groups in sufficiently high dimensions, i.e. in dimensions greater than the cohomological dimension of $G_F$.  This observation was formalized into the {\em Quillen-Lichtenbaum conjecture}.  This conjecture has attracted a great deal of attention.   Over the years, a number of special cases \cite{Suslin2},  \cite{HM} have been treated, and partial progress has been made \cite{DFST}.  Voevodsky in \cite{Voevodsky} proved results which led to the verification of the conjecture at $p=2$.  More recent   work of V. Voevodsky appears to have resolved this long standing conjecture at odd primes as well.   Moreover, it also appears to resolve the Beilinson-Lichtenbaum conjecture about the form of a spectral sequence which converges exactly to the  algebraic $K$-theory of $F$. 
 
 \noindent The existence of this spectral sequence does not, however, provide a homotopy theoretic model for the algebraic $K$-theory spectrum of the field $F$.  It is the goal of this paper to propose such a homotopy theoretic model, and to verify it in some cases.  We will construct a model which
 depends only on the complex representation theory (or twisted versions, when the roots of unity are not present in $F$) of $G_F$, without any other explicit arithmetic information about the field $F$,  and is therefore contravariantly functorial in $G_F$.  It is understood in algebraic topology that explicit space level models for spaces and spectra are generally preferable to strictly algebraic calculations of homotopy groups for the spaces.  Knowledge  of homotopy groups alone does not allow one to understand the behavior of various constructions and maps of spectra in the explicit way.  In this case, there are specific reasons why one would find such a model desirable.  
 
 \begin{itemize}
 \item{It appears likely that one could begin to understand {\em finite} descent problems, since the behavior of restriction and induction maps for subgroups of finite index is relatively well understood.}
 \item{It appears that the Milnor $K$-groups should  be identified with the homotopy groups of the so-called {\em derived completion} of the complex representation ring of $G_F$.   This would relate an explicitly arithmetic invariant (Milnor $K$-theory) with an explicitly group-theoretic invariant coming out of the representation theory.  The Bloch-Kato conjecture provides one such connection, relating Milnor $K$-theory with Galois cohomology.    These ideas provide  another one, arising out of the representation theory of the group  rather than from cohomology.}
 \item{It is well understood that representation theory of Galois groups plays a key role in many problems in number theory.  This construction provides another explicit link between representation theory and arithmetic in fields which it would be very interesting to explore. }
  \end{itemize}
  
  \noindent 

\noindent 
 The idea of the  construction is as follows.  We assume
for simplicity that $F$ contains an algebraically closed subfield.  A
more involved and more general version of the construction is discussed
in the body of the paper.  Let $F$ be a field, with separable closure
$\overline{F}$, and define the category of {\em descent data} for the
extension $F \subset \overline{F}$ (denoted by $V^G(\overline{F}))$ to
have objects the  finite dimensional $\overline{F}$-vector spaces $V$ with
$G_F$-action satisfying $\gamma (\overline{f}v) = \overline{f}^{\gamma}
\gamma (v)$, with equivariant $\overline{F}$-linear isomorphisms.    It
is standard descent theory  \cite{descent} that this category is equivalent to the
category of finite dimensional $F$-vector spaces.  On the other hand,
let  $Rep_F[G_F]$ denote the category of  finite dimensional continuous
representations of the profinite group
$G_F$.  There is a canonical homomorphism from $Rep_F[G_F]$ to
$V^G(\overline{F})$ obtained by applying $\overline{F} \column{\otimes}{F}
- $ and extending the action via the Galois action of $G_F$ on
$\overline{F}$.  We then have the composite

$$ A_F:KRep_k[G_F] \longrightarrow  KRep_F[G_F] \longrightarrow
K^G(\overline{F}) 
$$
 
 \noindent This map is of course very far from being an equivalence, since $\pi _0 KRep_k[G_F]$ is isomorphic to $R[G_F]$, a non-finitely generated abelian group, and $\pi _0 KF \cong \Bbb{Z}$.  However, we observe that  both spectra are {\em commutative $S$-algebras} in the sense of \cite{May et all}, and further that the map $A_F$ is a homomorphism of commutative $S$-algebras.  We describe in this paper a derived version of completion, which is applicable to any homomorphism of commutative $S$-algebras.  For such a homomorphism $f : A \rightarrow B$, we denote this derived completion by $A^{\wedge}_B$.  We also let $\Bbb{H}_p$ denote the mod-$p$ Eilenberg-MacLane spectrum.  $\Bbb{H}_p$ is also a commutatve $S$-algebra,  and we have  an evident commutative diagram of commutative $S$-algebras
 $$
 \begin{diagram}
 \node{KRep_k[G_F] } \arrow{e,t}{A_F} \arrow{s,t}{\varepsilon _p} \node{K^G \simeq KF} \arrow{s,t}{\varepsilon _p} \\
 \node{\Bbb{H}_p} \arrow{e,t}{id} \node{\Bbb{H}_p}
 \end{diagram}
 $$
 
 \noindent The vertical maps are induced by forgetful functors which forget the vector space structure and only retain the dimension mod $p$.  The completion construction is natural for such commutative squares, and we obtain a map of spectra $\alpha_F (p): KRep_k[G_F] ^{\wedge}_{{\Bbb{H}_p}} \longrightarrow KF^{\wedge}_{\Bbb{H}_p}$.  We will refer to $\alpha _F(p) $ as the ``representational assembly (at $p$)".  In this paper, we make this construction precise, and then study it, with the following results.  
 
\begin{itemize}
\item{ We prove that $\alpha _F(p) $ is an equivalence for fields of characteristic prime to $p$   containing a separably closed subfield with topologically finitely generated abelian  absolute Galois group.  This would include, for example, iterated power series fields $k((x_1))((x_2))\cdots ((x_n))$.  We regard this as a ``proof of concept", in that it shows that the map is highly non-trivial. For example, this result can be extended to prove that  for all fields of characteristic prime to $p$  containing a separably closed subfield, $\alpha _F(p)$ is surjective on $\pi _1$.    }
\item{We discuss how generally one would expect $\alpha _F(p) $ to be an equivalence.  Although we have no counterexamples, we suspect that for non-abelian absolute Galois groups, the statement above is inadequate and we would expect instead an equivariant form, involving the derived completions of Green functors, to be correct construction.  We discuss some supporting evidence that this should be the case.  }
\item{We discuss twisted versions, i.e. versions of $\alpha _F (p) $ which can be constructed when there is no separably closed subfield.  This can even include situations where the roots of unity are not in the ground field.  The case for field of positive characteristic is particularly nice, where there is a very  clean statement involving the $K$-theory of twisted group rings.  In this context, the corresponding conjecture for the case of a power series field over a finite field has been proved by G. Lyo in \cite{lyo}.  }
\item{We construct a possible geometric model for the derived representation theory, separate from the $k$-theory of the field.  For non-abelian profinite groups, the representation rings are complicated non-Noetherian rings, and evaluating derived completions is a daunting task.  For this reason, it appears sensible to approach the computation indirectly, by comparing it to a more geometric object.  One such object is {\em deformation $K$-theory}, a construction on discrete groups which builds a version of equivariant $K$-theory which takes into account the presence of continuous families of deformations of representations.  T. Lawson in \cite{lawson} has proved that for a finitely generated nilpotent group $\Gamma$, the derived representation theory of the pro-$p$ completion of $\Gamma$ is equivalent to the $p$-adic completion of its deformation $K$-theory.  This suggests strongly that approaches related to deformation $K$-theory may very well be useful models for the derived representation theory.  }

\end{itemize}
 
 \noindent In  future work, we plan to study various indirect methods for evaluation of the derived representation theory, including some which are particularly well suited for comparisons with the motivic methods.  We ultimately hope to be able to show, using these methods, that the homotopy groups of the derived completion at $p$  of the representation ring of the separable Galois group (regarded as a ring spectrum via the Eilenberg-MacLane construction) will be identified with $p$-completed Milnor $K$-theory, via an explicit homomorphism.  We also expect that these methods will increase our geometric understanding of this construction, which will take us in the direction of proving the conjecture as stated.  We also hope to formulate a conjecture which applies to all fields, even fields of characteristic zero where the roots of unity are not present in the ground field.

\noindent The ultimate hope is that the clarification of the relationship
between arithmetically defined descriptions of algebraic $K$-theory, such
as the motivic spectral sequence, with descriptions which involve the
Galois group and its representation theory directly, will shed more light
on arithmetic and algebraic geometric questions.  

\noindent The author wishes to express his thanks to a number of
mathematicians for useful conversations, including D. Bump, R. Cohen, W.
Dwyer, A. Elmendorf, E. Friedlander,  L. Hesselholt, M. Hopkins,  W.C.
Hsiang, J.F. Jardine, M. Levine, J. Li, I. Madsen, M.Mandell,  J.P. May, 
H. Miller, J. Morava, F. Morel,  K. Rubin, C. Schlichtkrull, V. Snaith, 
R. Thomason, R. Vakil, and C. Weibel.

\section{Preliminaries}
\label{Preliminaries}
We assume the reader to  be familiar with the category of spectra, as
developed in \cite{Maygeneral} and \cite{May et all} or \cite{Smith et
all}.  In either of these references, the category of spectra is shown to
possess a coherently commutative and associative monoidal structure,
called the smash product.  Of course, on the level of homotopy categories,
this monoidal structure is the usual smash product.  The presence of such
a monoidal structure makes it possible to define {\em ring spectra} as
monoid objects in
$\ul{Spectra}$.  It is also possible to define the notion of a
commutative ring spectrum, with appropriate higher homotopies encoding
the commutativity as well as associativity and distributivity.  In line
with the terminology of
\cite{May et all}, we will refer to these objects as $S$-algebras and
commutative
$S$-algebras, respectively. For the relationship with earlier notions of
$A_{\infty}$ and $E_{\infty}$ ring spectra, see \cite{May et all}. These same references will also introduce notions we will make use of, such as the Universal coefficient and the K\"{u}nneth spectral sequences.  See \cite{completion}, section 2.1,  for a discussion of the relevant notions as they will be used in the present paper. 
 
 \noindent We also  remind the reader of  the construction by Bousfield-Kan
of the $l$-completion of a space (simplicial set) at a prime $l$,
denoted by $X^{\wedge}_l$.  Bousfield and Kan construct a functorial cosimplicial
space
$T_l^{\cdot} X$, and define the $l$-completion of $X$ to be
$Tot(T_l^{\cdot})$.  This construction gives rise to a functorial
notion of completion at a prime $l$. It extends in a straightforward way to spectra by applying completion at $l$ levelwise. It is extended further  in \cite{completion} to the notion of completion of a module $M$ over a commutative $S$-algebra $A$ at a commutative $A$-algebra $B$, denoted $M^{\wedge} _B$.  The Bousfield-Kan construction is the special case of this construction where $A$ is the sphere spectrum and $B$ is the mod-$l$ Eilenberg-MacLane spectrum.  A number of  properties of this construction are established in \cite{completion}.  We will use them extensively.  We will use the special case of a group ring at one point in the paper.  The relevant results are the following. 
 
 \begin{GProposition} \label{grpequiv} Suppose that we have a homomorphism $\varphi:G. \rightarrow H.$ of abelian simplicial groups, so that $B\varphi$ induces an isomorphism on mod-$p$ homology groups. Then the evident homomorphism of commutative $S$-algebras 
$|A[G.]|^{\wedge}_{\Bbb{H}_p} \rightarrow |A[H.]|^{\wedge}_{\Bbb{H}_p}$ is a weak equivalence of $S$-algebras, where $\Bbb{H}_p$ is an algebra over the group ring via augmentation. 
\end{GProposition}
\begin{Proof} It follows immediately from \ref{grpringcomplete} below and the Atiyah-Hirzebruch spectral sequence that the map of cosimplicial spectra induced by $\varphi$ is a levelwise equivalence, since in this case the graded algebras 

$$ \pi _* (B  \column{\wedge}{A} B \column{\wedge}{A} \cdots \column{\wedge}{A} B)
$$
are $\Bbb{F}_p$-vector spaces.  It is standard that levelwise equivalences of cosimplicial spectra induce equivalences of total spectra. 
\end{Proof}
 
 \begin{GProposition}\label{grpringcomplete}
Suppose that we are given a commutative  $S$-algebra $A$ and an abelian simplicial group $G.$, and a homomorphism $A \rightarrow B$ of commutative $S$-algebras.  Suppose further that $i:G. \rightarrow E.$ is an inclusion of simplicial groups,  and that $E.$ is contractible.  Then the commutative $S$-algebra 
$$ \underbrace{B \column{\wedge}{|A[G.]|} B \column{\wedge}{|A[G.]|} \cdots \column{\wedge}{|A[G.]|} B}_{k \mbox{ factors }}
$$
 is equivalent to the commutative $S$-algebra

$$|\underbrace{ B \column{\wedge}{A} B \column{\wedge}{A} \cdots \column{\wedge}{A } B}_{k \mbox{ factors }} [E^k. / G^{k-1}.] |
$$
where $G^{k-1}. $ is included in  $E^k.$ via the homomorphism 
$(g_1, \ldots , g_{k-1}) \rightarrow (i(g_1), i(g_2)-i(g_1), i(g_3) - i(g_2 ), \ldots , i(g_{k-1})- i(g_{k-2}), i(g_{k-1}))$.   Note that the quotient group is a model for $BG^{k-1}.$.  
\end{GProposition}
\begin{Proof} Immediate from the definitions. 
\end{Proof}

We refer the reader to \cite{Quillen} for
results concerning $K$-theory spectra, notably the localization,
devissage, and reduction by resolution theorems.  These theorems apply
equally well to the completed versions $KA^{\wedge}_l$, where $l$ is a prime.  We will also
recall Suslin's theorem about the $K$-theory of an algebraically closed field \cite{Suslin}.

\begin{GTheorem} Let $k \rightarrow F$ be an inclusion of algebraically
closed fields of characteristic $p \neq l$ ($p$ may be 0). The the
natural map $Kk \rightarrow KF$ induces an equivalence $Kk^{\wedge}_l
\rightarrow KF^{\wedge}_l$. In fact, the proof shows that the map of
pro-spectra $Kk^{\wedge}_l \rightarrow KF^{\wedge}_l$ is a weak equivalence
in the sense that it induces an isomorphism of homotopy pro-groups.   
\end{GTheorem}

\noindent An important method for constructing spectra from combinatorial
data is via {\em infinite loop space machines} (see \cite{May1} or
\cite{Segal}), which are functors from the category of symmetric monoidal
categories to spectra.  The algebraic
$K$-theory functor is a prime example of this construction, since it can
obtained  by applying an infinite loop space machine to the symmetric
monoidal category of finitely generated projective modules over a ring. 
When the symmetric monoidal category has a coherently associative and
distributive second monoidal structure, such as tensor product of
modules, the spectrum constructed by an infinite loop space machine will
have the structure of an $S$-algebra.  If in addition the
second monoidal structure is coherently commutative, the spectrum will
be a commutative $S$-algebra. See \cite{elmendorf}
for these results.  Since the tensor product of finitely generated
projective modules over a commutative ring
$R$ is coherently commutative, we have 

\begin{GProposition} For any commutative ring $A$, the spectrum 
$KA$  is equipped with a commutative $S$-algebra structure  in a
canonical way.  The $S$-algebra structure is functorial for homomorphisms of commutative rings.  
\end{GProposition}

We also recall the ideas of equivariant stable homotopy theory. See for
example \cite{greenleesmay} or \cite{CarlssonG} for information about this
theory.  In summary, for a finite group $G$, there is a complete theory of
$G$-equivariant spectra which includes suspension maps for one-point
compactifications of all orthogonal 
representations of $G$.  The proper analogue of homotopy groups takes its
values in the abelian category of {\em Mackey functors}, which is a
suitably defined category of diagrams over a category whose objects are
finite $G$-sets, and whose morphisms include maps of $G$-sets, and also
transfer maps attached to orbit projections. The category of Mackey
functors admits a coherently commutative and associative tensor product,
which is denoted by $\Box$.  Consequently, by analogy with the theory of
rings, we define a {\em Green functor} to be a monoid object in the
category of Mackey functors.  The theory of ring functors has analogues
for most of the standard theorems and constructions of ring theory. In
particular, it is possible to define ideals, modules,  completions, tensor
products of modules over a Green functor $R$, and modules $Hom_R(M,N)$ for
any Green functor and modules $M$ and $N$ over $R$.  The category of
modules over a Green functor also has enough projectives, so
homological algebra can be carried out in this category.   It is also
easy to verify that by passing to direct limits, it is possible to
directly extend the ideas about Mackey and Green functors to profinite
groups.     See
\cite{Bouc} for background material about Mackey and Green functors.  

\noindent As mentioned above, the Mackey functor valued analogue of
homotopy groups plays the same role for equivariant spectra that ordinary
homotopy groups play for spectra.  For example, it is shown in
\cite{Lewis-May-McClure} that there exist Eilenberg-MacLane spectra
attached to every Mackey functor.  Moreover, a theory of ring and module
spectra in the equivariant category has been developed by May and Mandell
\cite{maymandell} and \cite{May4} in such a way that the Mackey-functor
valued homotopy group becomes a Green functor in an evident way, and that
the usual spectral sequences (K\"{u}nneth and Universal coefficient) hold
for modules in this category.


\section{Categories of descent data}

Let $A$ be any commutative ring, equipped with a group action by a profinite group
$G$.  We assume that the action of $G$ is continuous, in the sense that
the stabilizer of any element of $A$ is an open and closed subgroup of
finite index in $G$.  

\begin{gDefinition} By a {\em linear descent datum} for the pair $(G,A)$,
we will mean a finitely generated  $A$-vector space $M$, together with a
continuous action of $G$ on $M$ so that $g(am) = a^gg(m)$ for all $g \in
G$, $a\in A$, and $m  \in M$. We define two categories of linear descent
data, $V^G(A)$ and $V(G,A)$.  The objects of $V^G(A)$ are all linear
descent data for the pair $(G,A)$, and the morphisms are all equivariant
$A$-linear morphisms.  The objects of $V(G,A)$ are also all linear
descent data for $(G,A)$, but the morphisms are all $A$-linear morphisms
(without any equivariance requirements). The group $G$ acts continuously on the
category $V(G,A)$, by conjugation of maps (so the action is trivial on
objects), and the fixed point subcategory is clearly $V^G(A)$. Note that
both categories are symmetric monoidal categories under direct sum.  
\end{gDefinition}

\noindent We note that $\column{\otimes}{A}$ provides a coherently
associative and commutative monoidal structure on $V(G,A)$ and $V^G(A)$. 

\begin{gDefinition} We define the spectra $K^G(A)$ and $K(G,A)$ to be the
spectra obtained by applying an infinite loop space machine (\cite{May1}
or \cite{Segal})  to the symmetric monoidal categories of isomorphisms of 
$V^G(A)$ and
$V(G,A)$, respectively.  $K(G,A)$ is a spectrum with $G$-action, with
fixed point spectrum $K^G(A)$. The tensor product described above makes
each of these spectra into commutative $S$-algebras using the results of
\cite{elmendorf}. 
\end{gDefinition}

There are various functors relating these categories (and therefore
their $K$-theory spectra).  We have the fixed point functor
$$ (-)^G : V^G(A) \longrightarrow V^{\{e \}}(A^G) \cong A^G-mod
$$
defined on objects by $M \rightarrow M^G$. We also have the induction
functor

$$ A \column{\otimes}{A^G}-: A^G-mod \cong V^{\{e\}}(A^G) \longrightarrow 
V^G(A) 
$$
given on objects by $M \rightarrow A \column{\otimes}{A^G} M$. 
 The following
is a standard result in descent theory. See
\cite{descent} for details. 

\begin{gProposition} Let $F$ be a field, and suppose that we are given a continuous action of a profinite group $G$ on $F$, so that $F^G \hookrightarrow F$ is a Galois extension with
Galois group  $G$. Suppose further that $G$ is in fact the inverse limit of groups of order prime to the characteristic of $F$. Then  both
$(-)^G$ and
$F\column{\otimes}{F^G} -$ are equivalences of categories.  
\end{gProposition}
\begin{Proof} This is standard in the case when $G$ is finite, and follows by a direct limit argument in the profinite group case.  
\end{Proof}

\noindent We have an analogue of this result when $A$ is not a field.  We will write $V^G(A)_{proj}$ for the full subcategory of 
$V^G(A)$  on the modules which are finitely generated projective as modules over $A$.  

\begin{gProposition}\label{galoisring}
Suppose that $A$ is an infinite Galois extension of a ring $A_0 \subseteq A$, in the sense of \cite{demeyer}, and suppose that $A$ is a union of finite Galois extensions.  The group $\Gamma$  of automorphisms of $A$ over $A_0$ is a profinite group,  is equipped with a continuous action on $A$, and we have $A_0 = A^{\Gamma}$.  Suppose further that $A_0$ contains a field $k$, and that no subgroup of $\Gamma$ admits a continuous  homomorphism to the group $\Bbb{Z}/p\Bbb{Z}$, where $p$ is the characteristic of $k$.  Of course, this condition is vacuous if $char(k) =0$.  Then the restriction of the functor $(-)^{\Gamma}: V^{\Gamma}(A) \rightarrow A^{\Gamma}-mod$ to the subcategory 
$V^{\Gamma}(A)_{proj}$ factors through the category $Proj(A^{\Gamma})$ of finitely generated projective modules over $A^{\Gamma}$, and the restriction is an equivalence of categories. 
\end{gProposition}
\begin{Proof}  This is direct from Theorem 1.1, Ch. III of \cite{demeyer}, and follows from a direct limit argument in the general case. 
\end{Proof}

\noindent We also have 
\begin{gProposition} Let $F$ be a field. The category $V(G,F)$ is canonically equivalent to
the category $Vect(F)$ of finitely generated vector spaces over $F$. 
\end{gProposition}
\begin{Proof} The definition of the morphisms has no dependence on the
group action, and the result follows immediately. 
\end{Proof}

\noindent {\bf Remark:} Note that when the group action is trivial, i.e.
$G$ acts by the identity, and $A$ is a field, then the category $V^G(F)$ is just the category
of finite dimensional   continuous $F$-linear representations of $G$.  

\begin{gDefinition} In the case when the $G$-action is trivial, and $A$ is a field $F$,  we will
also write $Rep_F[G]$ for $V^G(F)$. 
\end{gDefinition}

\begin{gProposition} \label{repstructure} The functor

$$ Rep_{F^G}[G] \cong V^G(F^G) \stackrel{F \column{\otimes}{F^G}
-}{\longrightarrow} V^G(F) \cong V^{\{e \}}(F^G) \cong Vect(F^G)
$$
\noindent respects the tensor product structure, and $K(F^G)$ becomes an
algebra over the $S$-algebra $KRep_{F^G}[G]$.  
\end{gProposition}

\section{An example}

\noindent Let $k$ denote an algebraically closed field of characteristic
0, let $k[[x]]$ denote its ring of formal power series, and let $F = k((x))$
denote the field of fractions of $k[[x]]$.   In this
section, we wish to show that the representation theory of its absolute
Galois group can be used to create a model for the spectrum $K{k((x))}$. 

\noindent The field $F$ contains the subring
$A = k[[x]]$ of actual power series, and is obtained from it by inverting
$x$.  We begin by analyzing the spectrum
$KF$.  It follows from the localization sequence (see
\cite{Quillen}) that we have a fibration sequence of spectra

$$  Kk \longrightarrow KA \longrightarrow KF . 
$$

Further, the ring $A$ is a {\em Henselian local ring} with residue class
field $k$. It follows from a theorem of O. Gabber (see \cite{Gabber}) that
the map of spectra $KA \rightarrow Kk$ induces an isomorphism on
homotopy groups with finite coefficients, and consequently an equivalence
on
$l$-adic completions $\comp{KA}{l} \rightarrow \comp{Kk}{l}$. 
The fiber sequence above now becomes, up to homotopy equivalence,  a fiber
sequence

$$\comp{Kk}{l} \stackrel{i}{ \longrightarrow} \comp{Kk}{l} \longrightarrow
\comp{KF}{l}.
$$

\noindent All three spectra in this sequence become module spectra over
the commutative $S$-algebra $\comp{Kk}{l}$, and the sequence
consists of maps which are  $\comp{Kk}{l}$-module maps. Since the inclusion of rings $k \hookrightarrow F$ induces an isomorphism $\pi _0 \comp{Kk}{l} \cong \pi _0\comp{KF}{l} \cong
\Bbb{Z}_l$, we find that the  map $\pi _0 (i): \pi _0 \comp{Kk}{l}
\rightarrow \pi _0 \comp{KA}{l}$ is the zero map.  Since the module
$\pi _* \comp{Kk}{l}$ is cyclic over the ring $\pi _* \comp{Kk}{l}$, this
means that the inclusion  induces the zero map on all homotopy groups. 
The conclusion is that as a $\pi _*
\comp{Kk}{l}$-module,

 $$\pi _*\comp{KF}{l} \cong \pi _*\comp{Kk}{l} \oplus \pi
_*\Sigma \comp{Kk}{l}$$
\noindent where the second summand is topologicaly generated by the unit
$x$ viewed as an element of $K_1(F)$. In fact, the algebra structure is
also determined, since the square of the one dimensional generator is
zero.  Therefore the homotopy gropus are given by 

$$ \pi _*KF \cong \Lambda_{\Bbb{Z}_l[x]}(\xi)
$$

\noindent where $\Lambda$ denotes the Grassmann algebra functor,
where the polynomial generator
$x$ is in dimension 2, and where the exterior generator $\xi$ is in
dimension 1.

\noindent The following result is proved in \cite{Volklein}.  It is
derived from the fact that the algebraic closure of $k((x))$ is the field
of {\em Puiseux series}, i.e. the union of the fields
$k((x^{\frac{1}{n}}))$. 

\begin{gProposition}  The absolute Galois group $G$ of $k((x))$ is the
group
$\hat{\Bbb{Z}}$, the profinite completion of the group of integers.  
\end{gProposition}

\noindent Let $F$ denote the field $k((x))$, and let  $E$ denote
$\overline{k((x))}$. We have observed in
\ref{repstructure} that the spectrum
$\kay{F}$, which is equivalent to $K^{G}(E)$, becomes an algebra
spectrum over the $S$-algebra $KRep_{F}[G] \cong K^GF$.  We wish to
explore the nature of this algebra structure, and to use the derived
completion to demonstrate that the algebraic $K$-theory of $F$ can be
constructed directly from the representation theory of $G$ over an
algebraically closed field, e.g. $\Bbb{C}$.   

\noindent We have the map of $S$-algebras $KRep_F[G] \rightarrow K^G(E)
\cong \kay{F}$, induced by the functor 

$$  id_E \column{\otimes}{F} - : V^G(F) \longrightarrow V^G(E).
$$
We may compose this map with the canonical map $KRep_k[G] \rightarrow 
KRep_F[G]$ to obtain a map of $S$-algebras

$$ \hat{\alpha} :  KRep_k[G] \rightarrow K^G(E).
$$

\noindent We note that as it stands this map doesn't seem to carry much
structure.  

\begin{gProposition}  $\pi _* KRep_k[G] \cong R[G] \otimes ku_*$, where
$R[G]$ denotes the complex representation ring. (The complex
representation ring of a profinite group is defined to be the direct
limit of the representation rings of its finite quotients.) 
\end{gProposition}
\begin{Proof}  Since $k$ is algebraically closed of characteristic zero,
the representation theory of $G$ over $k$ is identical to that over
$\Bbb{C}$.  This shows that $\pi _0 KRep_k[G] \cong R[G]$.  In the
category
$Rep_k[G]$, every object has a unique decomposition into irreducibles,
each of which has $\Bbb{C}$ as its endomorphism ring.  The result follows
directly. 
\end{Proof}

\noindent For the rest of  this section, we will adopt the notational convention that all $K$-theory spectra and groups are $l$-adically completed, where $l$ is a prime.  So, when we write $KF$, it will denote the spectrum  $(KF)^{\wedge}_l$, and $K_*F$ will denote $\pi _* ((KF)^{\wedge}_l)$.  We will occasionally include the completion notation for emphasis, but not generally.

\noindent We also know from the above discussion that $K_* F \cong K_*k \oplus
K_{*-1} k $.  It is now easy to check that the map $KRep_k[G] \rightarrow
K^GE
\cong KF$ induces the composite

$$  R[G] \otimes K_*k \stackrel{\varepsilon \otimes id}{\longrightarrow}
K_*k \hookrightarrow K_*k \oplus K_{*-1}k \cong K_*F
$$

\noindent which does not appear to carry much information about
$\kay{F}$. However, we may use the  derived completion as follows.  As usual,
let $\Bbb{H}_l$ denote the mod-$l$ Eilenberg-MacLane spectrum.  We now
have a commutative diagram of $S$-algebras

$$
\begin{diagram}
\node{KRep_k[G]} \arrow{e,t}{\hat{\alpha}} \arrow{s,t}{\varepsilon}
\node{K^GE}
\arrow{s,b}{\varepsilon}
\\
\node{\Bbb{H}_l} \arrow{e,t}{id} \node{\Bbb{H}_l}
\end{diagram}
$$ 

\noindent where $\varepsilon$ in both cases denotes the augmentation map
which sends any vector space or representation to its dimension mod $l$.  The
naturality properties of the derived completion construction yield a
map 

$$  \alpha_F (l)  : KRep_k[G] ^{\wedge}_{\Bbb{H}_l} \rightarrow K^GE
^{\wedge}_{\Bbb{H}_l} \stackrel{\sim}{\rightarrow} KF
^{\wedge}_{\Bbb{H}_l}
$$ 
which we refer to as the {\em representational assembly map} for $F$.   The rightmost equivalence is induced by the functor which carries an object in $V^G(E)$ to its $G$-fixed point subspace.  The
goal of this section is to show that despite the fact that $\hat{\alpha}$
carries little information about $K_*F$, $\alpha_F (l) $ is an
equivalence of spectra.  

\begin{gProposition}\label{trivial} $KF^{\wedge}_{\Bbb{H}_l}$ is equivalent to 
the
$l$-adic completion of $\kay{F}$. 
\end{gProposition}
\begin{Proof} We note that the $KRep_k[G] $-module structure on $KF$ extends over the augmentation, i.e. there is a commutative diagram 
$$\begin{diagram} 
\node{KRep_k[G] \wedge KF }  \arrow{e} \arrow{s,b}{\varepsilon \wedge id} \node{KF} \arrow{s,t}{id} \\
\node{Kk \wedge KF} \arrow{e} \node{KF}  
\end{diagram}
$$
where the horizontal maps are the structure maps for the module structures.  Letting $A = KRep[G]$ and $B = Kk$, we write $KF_A$ and $KF_B$ for $KF$ regarded  as $A$ and $B$ module spectra, respectively.  We obtain a natural map 
$$  (KF _A) _{\Bbb{H}_l} \rightarrow (KF_B) _{\Bbb{H}_l}
$$
which is an equivalence by Proposition 3.2 (6) of \cite{completion}.  Next, let $KF_{S}$ denote $KF$ regarded as a module spectrum over the sphere spectrum $S^0$.  Then we similarly obtain a natural map $ (KF_S)_{\Bbb{H}_l} \rightarrow (KF_B)_{\Bbb{H}_l}$, which is an equivalence.  But, example 8.2 of \cite{completion} shows that $(KF_S)_{\Bbb{H}_l}$ is equivalent to the Bousfield-Kan $l$-adic completion of $KF$.
\end{Proof}

\noindent It remains to determine the structure of $\pi _* KRep_k[G]
^{\wedge} _{\Bbb{H}_l}$. 

\begin{gProposition}  There is an equivalence of spectra

$$KRep_k[G]^{\wedge}_{\Bbb{H}_l} \cong \comp{ku}{l} \wedge S^1_+.$$
In particular the homotopy groups are given by 
$\pi _i KRep_k[G]^{\wedge}_{\Bbb{H}_l} \cong \Bbb{Z}_l$ for all $i \geq
0$, and $\pi _i KRep_k[G]^{\wedge}_{\Bbb{H}_l}\cong 0 $ otherwise. 

\end{gProposition}
\begin{Proof} We first show that the $S$-algebra $KRep_k[G]$ may be identified with the group ring $Kk [\chi (\Bbb{Z}_l)] \cong Kk [ \Bbb{Z}/l^{\infty} \Bbb{Z}]$, where $\chi (\Bbb{Z}_l)$ denotes the group of complex characters of $\Bbb{Z}_l$, and $ \Bbb{Z}/l^{\infty} \Bbb{Z}$ denotes the union $\bigcup _n  \Bbb{Z}/l^{n} \Bbb{Z}$.  We observe that because $k$ is algebraically closed of characteristic zero, the irreducible representations are in bijective correspondence with the characters of $\Bbb{Z}_l$.  We let  $C[\chi (\Bbb{Z}_l)]$ denote the category whose objects are finite sets $X$ equipped with a function $\varphi _X : X \rightarrow \chi(\Bbb{Z}_l)$, and whose morphisms are set isomorphisms respecting the character valued functions.  $C[\chi (\Bbb{Z}_l)]$ admits a notion of sum (disjoint union) and of product (the product of objects $(X,\varphi _X) $ and $(Y, \varphi _Y)$ is the set $X \times Y$ equipped with the function $\varphi _{X \times Y}$ defined by $\varphi _{X \times Y}(x,y) = \varphi _X (x) \varphi _Y (y)$).  There is an evident functor from $C[\chi (\Bbb{Z}_l)]$ to $Rep_k[\Bbb{Z}_l]$, which carries disjoint unions to direct sums and products to tensor products, and which sends an object $(X, \varphi _X)$ in $C[\chi (\Bbb{Z}_l)]$ to the free $k$-vector space on the set $X$, with the representation on each element $x \in X$ specified by the character $\varphi _X (x)$.  Consequently, using the results of \cite{elmendorf},  we obtain a homomorphism of commutative $S$-algebras 
$$ S[\chi (\Bbb{Z}_l )] \longrightarrow KRep_k [\Bbb{Z}_l ].
$$
Since $KRep_k[ \Bbb{Z}_l] $ is an algebra over the commutative $S$-algebra $Kk$, we obtain a homomorphism $Kk [\chi (\Bbb{Z}_l )] \rightarrow KRep_k [\Bbb{Z}_l ]$, which is readily verified to be an equivalence of spectra.  

\noindent There is further a homomorphism $\chi (\Bbb{Z}_l) \rightarrow S^1.$, where  $S^1..$ denotes the singular complex of the circle group, and where the homomorphism is induced by the inclusion of topological groups
$$\chi (\Bbb{Z}_l) \cong \Bbb{Z}/l^{\infty} \Bbb{Z} \hookrightarrow  S^1.$$
\noindent as the group of $l$-power torsion points.  The  map of classifying spaces induces isomorphisms on mod-$l$ homology groups, and consequently by \ref{grpequiv} a weak  equivalence on derived completions.  \end{Proof}

\noindent We have now shown that the derived completion
$KRep_k[G]^{\wedge}_{\Bbb{H}_l}
$ has the same homotopy groups as the $K$-theory spectrum $\kay{F}$.  We
must now show that the representational assembly induces an isomorphism
on homotopy groups.  We  recall the localization sequence
\cite{Quillen} which allows us to compute $K_*(F)$. Consider the category
${\cal M} = Mod(k[[x]])$  of finitely generated  modules over
$k[[x]]$.  Let ${\cal T} \subseteq {\cal P}$ denote the full subcategory
of $x$-torsion modules.  Then the category $Vect(F)$ may be identified
with the quotient abelian category ${\cal M}/{\cal T}$, and we have the
usual localization sequence 

$$ \kay{{\cal T}} \rightarrow \kay{{\cal M}} \rightarrow \kay{{\cal
P}/{\cal T}}
$$
which is identified with the fiber sequence 

$$\kay{k} \rightarrow \kay{k[[x]]} \rightarrow \kay{F}.
$$

\noindent We now observe that there is a twisted version of this
construction.  Let ${\cal O}_E$ denote the integral closure of $k[[x]]$
in $E$.  ${\cal O}_E$ is closed under the action of the group $G$.  We
now define ${\cal E}^G$ to be the category whose objects are finitely
generated  ${\cal O}_E$ modules equipped with a $G$-action so
that $g(em) = e^g g(m)$.  Further, let ${\cal T}^G$ denote the full
subcategory of torsion modules.  

\begin{gProposition} The quotient abelian category ${\cal E}^G/{\cal
T}^G$ can be identified with the category of linear descent data
$V^G(E)$. 
\end{gProposition}
\begin{Proof} Straightforward along the lines of Swan's  proof  in \cite{Swan} that
${\cal M}/{\cal T} \cong Vect{F}$. \end{Proof}

\begin{gCorollary}\label{sequence} There is up to homotopy a fiber
sequence of spectra 

$$ K{\cal T}^G \longrightarrow K {\cal E}^G \longrightarrow K^G(E) \cong
KF.
$$
\end{gCorollary}

\noindent We note that for each of the categories ${\cal T}^G$, ${\cal
E}^G$, and $V^G(E)$, objects may be tensored with finite dimensional
$k$-linear representations of $G$ to give new objects in the category.  This construction has sufficient coherence so that the results of \cite{elmendorf} imply that 
 that each of the $K$-theory spectra are module spectra over
$KRep_k[G]$, and that the maps between them are morphisms of module spectra.

\begin{gLemma} The homotopy fiber sequence of \ref{sequence} is a
homotopy fiber sequence of $KRep_k[G]$-module spectra. 
\end{gLemma} 

\noindent  Derived completion at
$\Bbb{H}_l$ preserves fiber sequences of module spectra (\cite{completion}, Proposition 3.2 (3)), so we have a
homotopy fiber sequence 

$$(K{\cal T}^G)^{\wedge}_{\Bbb{H}_l} \longrightarrow ( K {\cal
E}^G)^{\wedge}_{\Bbb{H}_l}
\longrightarrow ( K^G(E))^{\wedge}_{\Bbb{H}_l}
\cong KF .
$$
The last equivalence is  \ref{trivial}. We will now show that 

\begin{itemize}
\item{$(K{\cal T}^G)^{\wedge}_{\Bbb{H}_l} \simeq *$}
\item{$ ( K {\cal
E}^G)^{\wedge}_{\Bbb{H}_l} \simeq KRep_k[G]^{\wedge}_{\Bbb{H}_l} . $}
\end{itemize}

\noindent The result will follow immediately.  

\begin{gProposition}\label{tortrivial}$(K{\cal
T}^G)^{\wedge}_{\Bbb{H}_l}
\simeq *. $
\end{gProposition}
\begin{Proof}  Let ${\cal T}^G_n $ denote  the
subcategory of finitely generated torsion $k[[t^{\frac{1}{n}}]]$-modules
$M$ equipped with $G$-action so that $g(fm) = f^g g(m)$ for $g \in G$, $f
\in k[[t^{\frac{1}{n}}]]$, and $m \in M$.  We have obvious functors

$$ \tau ^k_n = k[[t^{\frac{1}{k}}]] \column{\otimes}{k[[t^{\frac{1}{n}}]]}
- : {\cal T}_n^G \rightarrow {\cal T}_k ^G
$$

\noindent whenever $n$ divides $k$, and an equivalence of categories 

$$ \column{lim}{\rightarrow} {\cal T}_n^G \longrightarrow {\cal T}^G.
$$
A straightforward devissage argument now shows that the inclusion of
the full subcategory of objects of ${\cal T}_n^G$ on which
$x^{\frac{1}{n}}$ acts by $0$ induces an equivalence on $K$-theory
spectra. It follows that $K_*{\cal T}_n^G \cong R[G] \otimes K_*k$, and
moreover  that $K{\cal T}_n^G \simeq KRep_k[G]$ as module spectra over
$KRep_k[G]$. We will now need to analyze the map of $K$-theory spectra
induced by the functors $\tau _n^k$. Recall that $R[G] \cong
\Bbb{Z}[\Bbb{Q}/\Bbb{Z}]$. For all positive integers  $n$ and $s$, let
$\nu _{n,s}$ denote the element  in $R[G] $ corresponding to the element 
$ \displaystyle \mathop{\Sigma} _{i=0}^{s-1} [i/ns] \in
\Bbb{Z}[\Bbb{Q}/\Bbb{Z}]
$.  We claim that we have a commutative diagram

$$
\begin{diagram}
\node{K_* {\cal T}_n^G} \arrow{s,t}{\tau _n^{ns}} \arrow{e,t}{\sim}
\node{R[G]
\otimes K_*k}
\arrow{s,b}{(\cdot \nu _{n,s}) \otimes id} \\
\node{K_* {\cal T}_{ns}^G} \arrow{e,t}{\sim} \node{R[G] \otimes K_*k.}
\end{diagram}
$$ 

\noindent That we have such a diagram is reduced to a $\pi _0$
calculation, since the diagram  consists of
$R[G] \otimes K_*k$-modules which are obtained by extension of scalars from $R[G] \otimes K_0k$. It therefore suffices to check the
commuting of this diagram on the element $1 \in K_* {\cal T}_n^G$. But
$1$ under $\tau _n^{ns}$ clearly goes to the module $k[[x^{\frac{1}{ns}}]]
/ (x^{\frac{1}{n}})$. It is easily verified that as an element in the 
representation ring, this element is $\nu _{n,s}$. 

\noindent We are now ready to describe $\Bbb{H}_l \column{\wedge}{KRep[G]}
K{\cal T}^G$. Since smash products commute with filtering colimits,
we see that $$\Bbb{H}_l \column{\wedge}{KRep[G]}
K{\cal T}^G \simeq \mbox{ } \column{lim}{\rightarrow}\Bbb{H}_l
\column{\wedge}{KRep[G]} K{\cal T}_n^{G}.$$

\noindent But since $K_* {\cal T}_n^G$ is a free $R[G] \otimes
K_*k$-module of rank 1, 
$\Bbb{H}_l
\column{\wedge}{KRep[G]} K{\cal T}_n^{G} \cong \Bbb{H}_l$.  It is readily
checked that the induced map 

$$ id_{\Bbb{H}_l} \column{\wedge}{KRep[G]}(\cdot \nu _{ns}) \otimes id
$$
is multiplication by $s$ on $\pi _* \Bbb{H}_l$. It now clearly follows
that  $\pi _* \Bbb{H}_l \column{\wedge}{KRep[G]}
K{\cal T}^G = 0 $, and consequently that 

$$\pi _*\underbrace{ \Bbb{H}_l \column{\wedge}{KRep[G]}  \cdots \column{\wedge}{KRep[G]} \Bbb{H}_l}_{k \mbox{ factors}} \column{\wedge}{KRep[G]}
K{\cal T}^G = 0 
$$
\noindent for all $k$.  It now follows directly from the definition of the completion, together with the fact that the total spectrum of a cosimplicial spectrum which is levelwise contractible is itself contractible, that
$(K{\cal T}^G)^{\wedge}_{\Bbb{H}_l} \simeq *$.
\end{Proof}

\begin{gCorollary}\label{oneequiv} The map $(K{\cal E}^G)
^{\wedge}_{\Bbb{H}_l}
\rightarrow (K^G(E))^{\wedge}_{\Bbb{H}_l}$ is a weak equivalence of
spectra.
\end{gCorollary}
\begin{Proof}  Completion preserves fibration sequences of spectra. 
\end{Proof}

\noindent We next analyze $K_*{\cal E}^G$ in low degrees. Let
$\Bbb{N}^{\cdot}$ be the partially ordered set of positive integers,
where $m \leq n$ if and only if $m|n$.  We define a directed system of
$R[G]$-modules $\{A_n\}_{n>0}$ parametrized by $\Bbb{N}^{\cdot}$ by
setting 
$A_n = R[G]$, and where whenever $m|n$, we define the bonding map from
$A_m$ to $A_n$ to be multiplication by the element $ \nu
_{m,\frac{n}{m}}$.
We let $QR[G]$ denote the colimit of this module.  It follows from the
proof of \ref{tortrivial} that $K_*{\cal T}^G \cong QR[G] \otimes K_*k$. 
In particular, $K_*{\cal T}^G = 0 $ in odd degrees.  On the other hand,
we know that $K_* F \cong \Bbb{Z}_l$ in odd degrees, but that the elements in odd degree map injectively under  the connecting homomorphism in the
localization sequence. Consequently, we have 

\begin{gProposition}\label{odddegree} $K_* {\cal E}^G  =0$ in odd degrees.
In particular,
$K_1 {\cal E}^G = 0$. 
\end{gProposition}

\begin{gProposition}\label{evendegree} $K_0{\cal E}^G \cong R[G]$. 
\end{gProposition}
\begin{Proof} Objects of ${\cal E}^G$ are the same thing as modules over
the twisted group ring ${\cal O}_E\langle G \rangle$, and a devissage
argument shows that the inclusion 

$$ Proj({\cal O}_E\langle G \rangle ) \hookrightarrow Mod({\cal
O}_E\langle G \rangle ) $$
induces an isomorphism, so that we may prove the corresponding result for
finitely generated projective modules over ${\cal O}_E\langle G \rangle$.  Note that we have
a ring homomorphism 

$$ \pi : {\cal O}_E\langle G \rangle \longrightarrow k[G]
$$
\noindent given by sending all the elements $x^{\frac{1}{n}}$ to zero.
${\cal O}_E\langle G \rangle $ is complete in the $I$-adic topology,
where $I$ is the kernel of $\pi$.  
As in
\cite{Swan}, isomorphism classes of  projective modules over $ {\cal
O}_E\langle G \rangle $ are now in bijective correspondence via $\pi$
with the isomorphism classes of projective modules over $k[G]$, which
gives the result. 
\end{Proof} 

\begin{Remark}  \label{henselian} {\em There should be a better proof of this result, proceeding via a non-commutative version of the results of O. Gabber \cite{Gabber}. That is, one should prove directly that the reduction map ${\cal O}_E \langle G \rangle \rightarrow k[G]$ is an equivalence by proving that it is sufficiently ``Henselian".  It seems quite likely that such a proof could be made to work, which would mean that the proof would not require as input any computational information about the field in question.  }
\end{Remark} 

\begin{gCorollary}\label{twoequiv}  The functor ${\cal O}_E
\column{\otimes}{k} -  : Rep_k[G] \rightarrow {\cal E}^G$ induces a weak
equivalence on spectra.  Consequently, the map
$KRep_k[G]^{\wedge}_{\Bbb{H}_l}
\rightarrow  (K{\cal E}^G) ^{\wedge}_{\Bbb{H}_l}$ is a weak
equivalence. 
\end{gCorollary}
\begin{Proof} \ref{odddegree} and \ref{evendegree} identify the homotopy
groups of $K{\cal E}^G$ in degrees 0 and 1.  From the localization
sequence above, together with the fact that multiplication by the Bott
element induces isomorphisms $K_i {\cal T}^G \rightarrow K_{i+2} {\cal
T}^G$ and
$K_i F
\rightarrow K_{i+2} F$, it follows that it also induces an isomorphism 
$K_i {\cal E}^G \rightarrow K_{i+2} {\cal E}^G$. So, the map of the
statement of the corollary induces an isomorphism on homotopy groups,
which is the required result. 
\end{Proof}

\begin{gTheorem}  The map $\alpha _F (l) :
KRep_k[G]^{\wedge}_{\Bbb{H}_l} \rightarrow K^G(E)
^{\wedge}_{\Bbb{H}_l} \cong KF^{\wedge}_{\Bbb{H}_l}$ is a weak
equivalence of spectra. 

\end{gTheorem}
\begin{Proof} $\alpha _F (l)  $ is the composite of the maps of Corollary
\ref{oneequiv} and Corollary \ref{twoequiv}, both of which we have shown
are weak equivalences. 
\end{Proof}

\section{The case of an abelian absolute Galois group}\label{abeliancase}
We again adopt the convention that all $K$-theory spectra are assumed to be $l$-completed, and that $K$-groups are interpreted as homotopy groups of $l$-completed $K$-theory spectra.
We will now show how to extend the result of the previous section  to the case of geometric fields
(i.e. containing and algebraically closed subfield) $F$  of
characteristic zero with
$G_F$ a free pro-$l$ abelian group.  Let $F$ be geometric, and let $k
\subseteq F$ denote an algebraically closed subfield.  Let $A =
k[t_1^{\pm 1}, t_2^{\pm 1}, \ldots t_n ^{\pm 1}]$. We let $B$ denote the
ring obtained by adjoining all $p$-th power roots of unity of the
variables $t_i$ to $A$, so 
$$ B = \bigcup _s k[t_1^{\pm \frac{1}{l^s}},t_2^{\pm \frac{1}{l^s}},
\ldots  t_n^{\pm \frac{1}{l^{  s}}}]
$$
We choose a sequence of elements $\zeta _n$ in $k$ so that $\zeta _n$
is a primitive $l^n$-th root of unity, and so that $\zeta _n^l = \zeta
_{n-1}$. We define an action of the group $\Bbb{Z}_l^n$ on $B$ by $\tau
_i (t_i^{\frac{1}{l^n}}) = \zeta _n t_i^{\frac{1}{l^n}}$, and $\tau
_j(t_i^{\frac{1}{l^n}}) = t_i^{\frac{1}{l^n}}$ for $i \neq j$, where
$\{ \tau _1, \tau _2, \ldots , \tau _n \}$ is a set of topological
generators for $\Bbb{Z}_l^n$.

\begin{gProposition} The $l$-completed $K$-theory groups of $A$ and $B$
are given by 

\begin{itemize}
\item{$K_*A \cong \Lambda_{K_*k}(\theta_1, \theta _2, \ldots, \theta _n)$}
\item{$K_*B \cong K_*k$}
\end{itemize}
\end{gProposition}
\begin{Proof} The first result is a direct consequence of the formula for
the $K$-groups of a Laurent polynomial ring.  The second also follows
from this fact, together with  analysis of the behavior of the $l$-th
power map on  a Laurent  extension.\end{Proof}

\noindent We now may construct a representational assembly in this case
as well, to obtain a representational assembly 

$$  \alpha _{L} (l) :KRep_k[G]^{\wedge}_{\Bbb{H}_l}
\longrightarrow K^G(B)^{\wedge}_{\Bbb{H}_l} \cong KA ^{\wedge}_{\Bbb{H}_l} \cong KA^{\wedge}_l
$$
where the first equivalence follows from \ref{galoisring}, and the second in an identical fashion to that in the proof of \ref{trivial}.

\begin{gProposition} $\alpha_{L} (l)$ is an equivalence of
spectra.
\end{gProposition} 
\begin{Proof} We first consider the case $n = 1$, and the field $F =
k((x))$, where we have already done the analysis.  It is immediate that
the inclusion $A \rightarrow F$ defined by $t \rightarrow x$ extends to
an equivariant homomorphism $B \rightarrow E$ of $k$-algebras, and
consequently that we get a commutative diagram

$$
\begin{diagram} \node{KRep_k[G]} \arrow{se,t}{\alpha _F (l) }
\arrow{s,t}{\alpha_{L}(l)}
\node{}\\
\node{KA} \arrow{e,t}{\cong} \node{KF}
\end{diagram}
$$
It is readily checked that the inclusion $KA \rightarrow KF$ is a weak
equivalence, and so that the horizontal arrow in the diagram is a weak
equivalence.  On the other hand, we have already shown that $\alpha
_{F}(l)$ is a weak equivalence of spectra, which shows that 
$\alpha_{L}(l)$ is an equivalence.  For larger values of $n$,
the result follows by smashing copies of this example together over the
coefficient $S$-algebra $K_*k$. 
\end{Proof}
 
\begin{Remark} {\em As in Remark \ref{henselian}, we believe that there is another proof of this result, this time following from the homotopy property of \cite{Quillen}.  The idea would be to observe that the group action preserves the polynomial subrings $k[t_1^{\frac{1}{l^s}}, \ldots , t_n^{\frac{1}{l^s}}]$, and that a  corresponding twisted polynomial ring should have $K$-theory computable using Quillen's results.}  \end{Remark}

\noindent  Now consider any geometric field $F$ of
characteristic zero, with $G_F \cong \Bbb{Z}_l^n$, and let $k$ be an
algebraically closed subfield. It is a  direct consequence of Kummer
theory that there is a family of elements $\{ \alpha _1, \alpha _2,
\ldots , \alpha _n \}$ so that 
$E =
\overline{F}$ can be written as 

$$  E = \bigcup _n F( \sqrt[l^n]{\alpha _1}, \sqrt[l^n]{\alpha _2},
\ldots, \sqrt[l^n]{\alpha _n} )
$$ 

\noindent It is therefore further clear that we may define a 
$k$-algebra 
homomorphism $i:B \rightarrow E$ by setting $i(t_i^{\frac{1}{l^n}}) =
\sqrt[l^n]{\alpha _i}$, where $\sqrt[l^n]{\alpha _i}$ is chosen as the
choice of $l^n$-th root so that the action of $\tau _i$ is via
multiplication by $\zeta _n$. It follows that we obtain a commutative
diagram

$$
\begin{diagram} \node{KRep_k[G]^{\wedge}_{\Bbb{H}_l}}
\arrow{s,t}{\alpha_{L}(l)} \arrow{se,t}{\alpha _F(l)} \\
\node{KA} \arrow{e,t}{Ki} \node{KF}
\end{diagram}
$$
Since we now know that $\alpha_{L}(l)$ is a weak equivalence, it
will now suffice to show that $Ki$ is a weak equivalence of spectra.  But the work of Geisser-Levine \cite{Geisser} shows that if the Bloch-Kato conjecture holds, then the $l$-completed version of the Bloch-Lichtenbaum spectral sequence converging to the algebraic $K$-theory has as its $E_2$-term an exterior algebra on generators in degree one, corresponding to the elements $\alpha _1, \alpha _2, \ldots , \alpha _n$, from which this result follows directly.     

\begin{gProposition} If the Bloch-Kato conjecture holds, then $\alpha _F(l)$ is an equivalence of spectra for
geometric fields of characteristic prime to $l$  with finitely generated abelian separable  Galois
group. 
\end{gProposition}
\section{Surjectivity on $\pi _1$ in the geometric case}

\noindent Let $F$ be any field containing an algebraically closed subfield $k$, and  suppose $l$ is prime to the characteristic of $k$. We claim that the map $\alpha _F (l)$ induces a surjection on $\pi _1$.  

\begin{Proposition}
For any field $F$, we have $$\pi _1((KF)^{\wedge}_l) \cong (F^{*})^{\wedge}_l$$
where the right hand term denotes the usual algebraic $l$-adic completion of the group $F^*$.  
\end{Proposition}
\begin{Proof}  Since $\pi _0 (KF) \cong \Bbb{Z}$ is a finitely generated abelian group, this result  follows from Proposition 5.1 in  Ch VI of \cite{Bousfield}. 
\end{Proof} 

\noindent Functoriality of the assembly construction now  gives us a density result.  
\begin{Proposition} Suppose $F$ contains an algebraically closed subfield $k$, and that $l$ is prime to $char(k)$.  Then the image of the homomorphism
$$\pi _1 (\alpha _F (l)):  \pi _1 ( KRep_k[G] ^{\wedge}_{\Bbb{H}_l}) \rightarrow \pi _1 ((KF)^{\wedge}_l ) \cong (F^*)^{\wedge}_l
$$
is dense in the profinite topology on $(F^*)^{\wedge}_l$. In 
\end{Proposition} 
\begin{Proof}  It clearly suffices to prove that for any $n$ and any element 
$x \in F^* /l^nF^*$, there is an element $ \xi \in \pi _1 ( KRep_k[G] ^{\wedge}_{\Bbb{H}_l})$ so that $\rho(\pi _1 (\alpha _F (l))(\xi) ) = x$, where $\rho $ is the projection $F^*  \rightarrow F^* /l^nF^*$.  In order to do this, we select $\overline{x} \in F^*$ so that $\rho (\overline{x}) = x$. We next  consider the extension $\bigcup _n F(x^{\frac{1}{l^n}}$ of $F$.  This extension is preserved under the action of the separable Galois group $G$, and therefore yields a particular topologically cyclic quotient of $G$, which we denote by $Q$.   We note that $Q$ is also a pro-$l$ group, and is therefore isomorphic to $\Bbb{Z}_l$.  Select a particular identification  $\theta: Q \simeq \Bbb{Z}_l$.   Next,  define a ring homomorphism $A = k[t^{\pm 1}] \rightarrow F$ by sending the variable $t$ to $\overline{x}$.  Choosing a sequence of primitive $l^n$-th roots of unity $\zeta _n$ as in the  beginning of Section \ref{abeliancase}, we extend this map to a ring homomorphism 
$$ i: B =  \bigcup _n k[t^{\pm \frac{1}{l^n}}] \rightarrow F
$$
via the requirement that $t^{\pm \frac{1}{l^n}}$ be sent to the unique $l^n$-th root of $x$ on which the topological generator $\theta ^{-1}(1)$ acts by multipication by $\zeta _n$.  It is now clear from the construction that the diagram 

$$    
\begin{diagram}
\node{KRep_k[\Bbb{Z}_l]^{\wedge}_{\Bbb{H}_l}} \arrow{e,t} {\alpha _L(l)} \arrow{se,b}{\alpha _F(l)} \node{KA^{\wedge}_l} \arrow{s} \\
\node{}  \node{KF^{\wedge}_{l}}
\end{diagram}
$$
where the vertical arrow is induced by the ring homomorphism $A \rightarrow F$ defined above.  It was proved in Section \ref{abeliancase} that $\alpha _L(l)$ is an equivalence, hence induces an isomorphism on $\pi _1$.  The result now follows immediately.  
\end{Proof}

\noindent The actual surjectivity now follows from 

\begin{Proposition}  The groups are $\pi _t ( KRep_k[G]^{\wedge}_{\Bbb{H}_l})$ are all $l$-complete.  
\end{Proposition}
\begin{Proof}  It follows from the definition of the derived completions that the derived completion is the total spectrum of a cosimplicial spectrum which in every level is a mod-$l$  Eilenberg-Maclane spectrum, hence is $l$-complete in the sense of \cite{Bousfield}.  It is standard that the total space of such a cosimplicial spectrum itself has $l$-complete homotopy groups.  
\end{Proof}

\section{Representational assembly  in the geometric case}
\label{secthree}

Suppose now that we are in the case of a {\em geometric field}, i.e. a
field $F$  containing an algebraically closed subfield $k$.  Let $E$
denote the algebraic closure of $F$, and let $G$ denote the absolute
Galois group. As in the case of the example of the previous section, we
have a composite functor 

$$Rep_k[G] \longrightarrow Rep_F[G] \cong V^G(F) \longrightarrow V^G(E)
$$
and the induced map of completed  spectra

$$ KRep_k[G]^{\wedge}_{\Bbb{H}_l} \longrightarrow
K^G(E)^{\wedge}_{\Bbb{H}_l} \cong KF^{\wedge}_{\Bbb{H}_l}
$$

\noindent We may ask whether this map is an equivalence as occurred in
the case of the example in the previous section.  In order to discuss the
plausibility of this kind of statement, let us examine what occurs in the
example, when $G = \Bbb{Z}_l$.  In this case, $\pi _0 KRep_k[G] \cong
\Bbb{Z} [\Bbb{Q} /\Bbb{Z}]$.  This is a very large, non-Noetherian
ring.  However, after derived completion, we find that $\pi _0
KRep_k[G]^{\wedge}_{\Bbb{H}_l} \cong \Bbb{Z}_l$.  On the other hand, we
have from Theorem 7.1 of \cite{completion} that $\pi _0
KRep_k[G]^{\wedge}_{\Bbb{H}_l}$ is isomorphic to the ordinary (not
derived) completion of the ring $R[G]$ at the ideal $J = (l) + I$,
where $I$ denotes the  augmentation ideal in $R[G]$.  So the conclusion
is that
$
\column{lim}{\leftarrow}R[G]/J^n
\cong \Bbb{Z}_l \cong (R[G]/I) ^{\wedge}_l$.  The reason this happens is
the following result about the ideal $I$. 

\begin{gProposition}\label{div}  We have that $I/I^2$ is a divisible
group, and that
$I^k = I^{k+1}$ for all $k$. 
\end{gProposition}
\begin{Proof}  It is a standard result for any group $G$ that in the
group ring $\Bbb{Z}[G]$, the augmentation ideal $I[G]$ has the property
that
$I[G]/I[G]^2 \cong G^{ab}$.  Since $\Bbb{Q}/\Bbb{Z}$ is abelian, we find
that 
$I/I^2 \cong \Bbb{Q}/\Bbb{Z}$, which is divisible.  For the second
statement, we note that $I^k/I^{k+1}$ is a surjective image of 
$\underbrace{I/I^2 \otimes \ldots \otimes I/I^2}_{k \mbox{ factors}}$. 
But it is clear that $\Bbb{Q}/\Bbb{Z} \otimes \Bbb{Q} / \Bbb{Z} = 0$, so
$I^k/I^{k+1} = 0$. This gives the result.  
\end{Proof} 

\begin{gCorollary}\label{newdiv} We have that $((l)+I)^k = (l^k) + I$.
Equivalently, for any element $\theta \in I$, and any positive integers
$s$ and $t$, there is an $\eta \in I$ and $\mu \in I^k$ so that 

$$   \theta = l^s \eta + \mu . 
$$ 
\end{gCorollary}

\begin{Proof} We have 
$$ ((l)+I)^k = \mathop{\Sigma} _t l^t I^{k-l} = (l^k)+l^{k-1}I +I^2. 
$$
But \ref{div} implies that $l^{k-1}I + I^2 = I$, which gives the result.  
\end{Proof}

\noindent It now follows that $R[G]/J^k \cong R[G]/(I+(l^k)) \cong
\Bbb{Z}/l^k \Bbb{Z}$, and therefore that the completion of $R[G]$ at
the ideal $J$ is just $\Bbb{Z}_l$.  

\noindent This result extends to torsion free abelian profinite groups. 
That is, any such group also has the property that if $I$ is the
augmentation ideal in the representation ring, then
$I^k = I^2$ for
$k
\geq 2$, and that $I/I^2$ is divisible.  However, for non-abelian Galois
groups, there appears to be no obvious reason why such a result should
hold, and indeed we believe that it does not.  However, there is a
modification of this statement which is true for absolute Galois groups
and which suggests a conjecture which can be formulated for any geometric
field. 

\noindent  Recall that the representation ring of a finite group is not
just a ring, but is actually a part of a {\em Green functor}.  This was
discussed in section \ref{Preliminaries}, where there is a discussion of
Mackey and Green functors. See \cite{Bouc} for a thorough discussion of
these objects.   A Mackey functor for a group
$G$ is a functor from a category of orbits to abelian groups.  In the
case of the representation ring, there is a functor ${\cal R}$  given on
orbits by
${\cal R}(G/K) = R[K]$.  The maps induced by projections of orbits induce
restriction maps on representation rings, and transfers induce
inductions. This functor is actually a commutative {\em Green functor},
in the sense that there is a multiplication map ${\cal R} \Box {\cal R }
\rightarrow {\cal R}$, which is associative and commutative. Also for any
finite group $G$, there is another Green functor ${\cal Z}$ given on
objects  by ${\cal Z}(G/K) =
\Bbb{Z}$, and for which projections of orbits induce the identity map and
where transfers associated to projections  induce multiplication by the
degree of the projection.  We may also consider ${\cal Z} /l {\cal Z}$,
which is obtained by composing the functor ${\cal Z}$ with the
projection $\Bbb{Z} \rightarrow  \Bbb{Z}/l\Bbb{Z}$. The {\em
augmentation} is the morphism of Green functors ${\cal R} \rightarrow
{\cal Z}$ which is given on an object $G/K$ by the augmentation of
$R[K]$. The mod-$l$ augmentation $\varepsilon$  is the composite 
${\cal R} \rightarrow {\cal Z} \rightarrow {\cal Z}/l{\cal Z}$. The
theory of Green functors is entirely parallel with the theory of rings,
modules, and ideals.  We may therefore speak of the {\em
 augmentation ideal} ${\cal I}$ and the ideal ${\cal J} = (l) + {\cal
I}$, as well as the powers of these ideals. We can therefore also speak
of completion at a Green functor ideal. We also observe that the theory
of Mackey and Green functors extends in an obvious way to profinite
groups, by considering the category of finite $G$-orbits.   We want to
identify a class of profinite groups
$G$ for which 

$${\cal R}^{\wedge}_l = \column{lim}{\leftarrow} {\cal R} / {\cal J}^n
$$

\noindent is isomorphic to ${\cal Z}^{\wedge}_l = \Bbb{Z}_l \otimes 
{\cal Z}$. 

\begin{gDefinition} Let $G$ be a profinite group.  We say $G$ is {\em
totally torsion free} if every subgroup of finite index has a torsion
free abelianization.   
\end{gDefinition}

{\bf Example:} Free  profinite groups and free profinite $l$-groups
are totally torsion free. 

{\bf Example:} Free profinite abelian and profinite $l$-abelian groups. 

{\bf Example:} Let $\Gamma$ denote the integral Heisenberg group, i.e.
the group of upper triangular integer matrices with ones along the
diagonal.  Then the profinite and pro-$l$ completion of $\Gamma$ is not
totally torsion free.

\begin{gProposition} Let $G$ be the absolute Galois group of a
geometric field
$F$. Let
$G_l$ denote the maximal pro-$l$ quotient of $G$.  Then $G_l$ is totally
torsion free. 
\end{gProposition}
\begin{Proof}    Consider any subgroup $K$ of finite
 index in $G_l$.  Then let $E$ denote the extension of $F$
corresponding to $K$.  Then the abelianization of $K$ corresponds to
the maximal pro-$l$ Abelian extension of $E$. Let $\Bbb{N}^{\cdot}$
denote the partially ordered set of positive integers.  Define a functor $\Phi$ from
$\Bbb{N}^{\cdot}$ to abelian groups by $\Phi (n) = E^*/(E^*)^{l^n}$, and
on morphisms by $\Phi (m \leq  m+n) = E^*/(E^*)^{l^m}
\stackrel{(-)^{l^n}}{\rightarrow} E^*/(E^*)^{l^{m+n}}$. Let ${\cal E}$ denote
the direct limit over $\Bbb{N}^{\cdot}$ of $\Phi$.  Kummer theory then
asserts the existence of a perfect pairing 

$$ K^{ab} \times {\cal E} \rightarrow \Bbb{Q}/\Bbb{Z}
$$

\noindent So $K^{ab} \cong Hom({\cal E}, \Bbb{Q}/\Bbb{Z})$.  The group
${\cal E}$ is clearly $l$-divisible, and it is easily verified that the
$\Bbb{Q}/\Bbb{Z}$-dual of an $l$-divisible group is $l$-torsion free.   
\end{Proof}

\noindent We have the following result concerning totally torsion free
profinite groups.

\begin{gProposition} \label{torfreerep} Let $G$ be a totally torsion free
$l$-profinite group.  Then the natural map ${\cal R}^{\wedge}_l
\rightarrow {\cal Z}^{\wedge}_l$ is an isomorphism. Consequently, for $G$
the maximal pro-$l$ quotient of the absolute Galois group of a geometric
field, we have that ${\cal R}^{\wedge}_l \cong {\cal Z}^{\wedge}_l$. 
\end{gProposition}

\begin{Proof}  We will verify that ${\cal R}^{\wedge}_l(G/G) \rightarrow
{\cal Z}^{\wedge}_l(G/G)$ is an isomorphism.  The result at any $G/K$
for any finite index subgroup will follow by using that result for the
totally torsion free group $K$. It will suffice to show that for every
finite dimensional representation
$\rho$ of $G/N$, where
$N$ is a normal subgroup of finite index, and every choice of positive
integers
$s$ and
$t$, there are elements $x \in {\cal R}(G/G)$ and $y \in {\cal I}^t(G
/G)$ so that $[dim\rho] - [\rho] = l^s x + y$.  We recall {\em
Blichfeldt's theorem} \cite{Serre2}, which asserts that there is a
subgroup
$L$ of $G/N$ and a one-dimensional representation $\rho _L$ of $L$ so that
$\rho$ is isomorphic to the representation of $G/N$ induced up from $\rho
_L$. It follows that $[dim\rho] - [\rho] = i_L^{G/N}(1-\rho _L)$. Let
$\overline{L}$ denote the subgroup $\pi ^{-1}L \subseteq G$, where $\pi
: G \rightarrow G/N$ is the projection, and let $\rho _{\overline{L}} =
\rho_L \circ \pi$.   Then we clearly also have
$[dim\rho] - [\rho] = i_{\overline{L}}^{G}(1-\rho _{\overline{L}})$. 
Now, $1 - \rho _{\overline{L}}$ is in the image of $R[\overline{L}^{ab}]
\rightarrow R[\overline{L}]$, and let the corresponding
one-dimensional representation  of $\overline{L}^{ab}$ be $\rho
_{\overline{L}^{ab}}$. Since $\overline{L}^{ab}$ is abelian and torsion
free (by the totally torsion free hypothesis), we may write 
$1-\rho
_{\overline{L}^{ab}} = l^s \xi + \eta$, where $\eta \in
I^t(\overline{L}^{ab})$ and where $\xi \in R[\overline{L}^{ab}]$, by
\ref{newdiv}. This means that we may pull $\xi $ and $\eta$ back along the
homomorphism 
$\overline{L} \rightarrow \overline{L}^{ab}$, to get elements
$\overline{\xi} \in R[\overline{L}] $ and $\overline{\eta} \in
I^t(\overline{L})$ so that $\rho _{\overline{L}} = l^s \overline{\xi} +
\overline{\eta}$. Since ${\cal I}^t$ is closed under induction, we have
that $i_{\overline{L}}^G(\overline{\eta}) \in {\cal I}^t(G/G)$. Now we
have that $[dim \rho] - \rho = i_{\overline{L}}^G (1 - \rho
_{\overline{L}}) = p^s i_{\overline{L}}^G \overline{\xi} +
i_{\overline{L}}^G \eta$.  The result follows.   
\end{Proof}

We recall the relationship between Mackey functor theory and equivariant
stable homotopy theory.  The natural analogue for homotopy groups in the
world of equivariant spectra takes its values in the category of Mackey
functors.  The 
Adams spectral sequence and the algebraic to geometric spectral sequence
have $E_2$-terms which are computed using derived functors of Hom and
$\Box$-product, which we will denote by {\em Ext} and {\em Tor} as in the
non-equivariant case.  We will now observe that the constructions we have
discussed are actually part of equivariant spectra.  

\begin{gProposition}  For any profinite group $G$, there is a stable
category of $G$-spectra, which has all the important properties which the
stable homotopy theory of $G$-spectra for finite groups has.  In
particular, there are fixed point subspectra for every subgroup of finite
index as well as transfers for finite $G$-coverings.  The ``tom Dieck''
filtration holds as for finite groups.  The homotopy groups of a
$G$-spectrum form a Mackey functor.  Moreover, the homotopy groups of a
$G$ $S$-algebra are a Green functor, and the homotopy groups of a
$G$-module spectrum are a module over this Green functor.  
\end{gProposition}

\begin{Proof}  For any homomorphism $f : G \rightarrow H$ of finite
groups, there is a  pullback functor $f^*$ from the category of
$H$-spectra to the category of $G$-spectra.  Hence, for a profinite group
$g$, we get a direct limit system of categories parametrized by the
normal subgroups of finite index, with the value at $N$ being the
category of $G/N$-spectra.  A $G$ spectrum can be defined as a family of
spectra ${\cal S}_{N}$ in the category of $G/N$-spectra, together with
isomorphisms $$(G/N_1 \rightarrow G/N_2)^*({\cal S}_{N_2})
\stackrel{\sim}{\longrightarrow } {\cal S}_{N_1}. $$  This general
construction can be made into a theory of $G$ spectra with the desired
properties. 
\end{Proof}

\begin{gProposition} Let $F$ be a field, with $E$ its algebraic closure,
and $G$ the absolute Galois group.  Then there is a $G$-$S$-algebra with
total spectrum $K(G,E)$, with fixed point spectra $K(G,E)^H \cong
K^H(E)$.  The attached Green functor is given by 

$$ G/L \rightarrow \pi _*K(G,E)^L \cong \pi _* K^L(E) \cong K_*(E^L). 
$$  

\noindent Similarly, there is a $G$-$S$-algebra with total spectrum
$K(G,F)$, and with fixed point spectra $K(G,F)^L \cong K^L(F) \cong
KRep_F[L]$. (Note that the $G$-action on $F$ is trivial.) In this case,
the associated Green functor is given by 

$$G/L \rightarrow \pi _* K(G,F)^L \cong \pi _* K^L(F) \cong K_* Rep_F[L]
$$

\noindent In the case when $F$ contains all the $l$-th power roots of
unity, we find that the Mackey functor attached to $K(G,F)$ is $ku_*
\otimes {\cal R}$.  The functor
$E
\column{\otimes}{F} -$ induces a map of
$G$-$S$-algebra $K(G,F) \rightarrow K(G,E)$, which induces the
functors of Proposition \ref{repstructure} on fixed point spectra.

\end{gProposition}

\begin{Proof} This result is essentially a consequence of the equivariant
infinite loop space recognition principle \cite{Shimakawa}. 
\end{Proof}

\noindent We also recall the results of \cite{Lewis-May-McClure}, where it
was shown that in the category of $G$-equivariant spectra, there is an 
Eilenberg-MacLane spectrum for every Mackey functor.  We will let ${\cal
H}_l$ denote the Eilenberg-MacLane spectrum attached to the Green functor
${\cal Z}/l{\cal Z}$.  Derived completions of homomorphisms of
$S$-algebras are defined as in the non-equivariant case. 

\begin{gProposition} Let $F$ be a geometric field, with the
algebraically closed subfield $k$.  There is a commutative
diagram of
$G$-$S$-algebras

$$
\begin{diagram}
\node{K(G,k)} \arrow{e,t}{E\column{\otimes}{k}} \arrow{se,t}{\varepsilon}
\node{K(G,E)}
\arrow{s,b}{\varepsilon}\\
\node{} \node{{\cal H}_l}
\end{diagram}
$$ 

\noindent where  the maps $\varepsilon$ are given by ``dimension mod
$l$''.  Consequently, there is a map of derived completions

$$ \alpha^{rep}: K(G,k) ^{\wedge}_{{\cal H}_l} \longrightarrow K(G,E)
^{\wedge}_{{\cal H}_l}
$$

\end{gProposition}

It follows from the equivariant algebraic to geometric spectral sequence
that $\pi _0 K(G,k) ^{\wedge}_{{\cal H}_l} \cong {\cal R}^{\wedge}_l$,
and  by Proposition \ref{torfreerep} this is isomorphic to ${\cal
Z}^{\wedge}_l$, which is in turn isomorphic to the Green functor $\pi _0
K(G,E)$.  This makes  plausible  the following
conjecture.

\begin{gConjecture} \label{assemblyconj} Let $F$ be a geometric field. The
representational assembly map 

$$  \alpha ^{rep}: K(G,k)^{\wedge}_{{\cal H}_l} \rightarrow K(G,E)
^{\wedge}_{{\cal H}_l}
$$
is an equivalence of $G$-$S$-algebras. On $G$-fixed point spectra, we
have an equivalence of $S$-algebras

$$ (\alpha ^{rep})^G : (K(G,k)^{\wedge}_{{\cal H}_l})^G \longrightarrow
(K(G,E)
^{\wedge}_{{\cal H}_l})^G \cong K(E^G) \cong KF.  
$$
\end{gConjecture}

\section{Representational assembly in the twisted case}

In the case of a field $F$ which  does not contain an algebraically closed
subfield, it is not as easy to construct a model for the $K$-theory
spectrum of the field. As usual, let $F$ denote a field, $E$ its
algebraic closure, and $G$ its absolute Galois group.  We observe first
that for any subfield of $F^{\prime} \subseteq E$ which is closed under
the action of $G$, we obtain as above  a diagram

$$
\begin{diagram}
\node{K(G,F^{\prime})} \arrow{e,t}{E\column{\otimes}{F^{\prime} }} \arrow{se,t}{\varepsilon}
\node{K(G,E)}
\arrow{s,b}{\varepsilon}\\
\node{} \node{{\cal H}_l}
\end{diagram}
$$ 
 and consequently a map 

$$  K(G,F^{\prime}) ^{\wedge}_{{\cal H}_l} \longrightarrow K(G,E)
^{\wedge}_{{\cal H}_l}. 
$$

\noindent We may conjecture (as an extension of Conjecture
\ref{assemblyconj} above) that this map is always an equivalence, and
indeed we believe this to be true.  In the case when $F^{\prime}$ is
algebraically closed, we  believe  that the domain of the map has what
we would regard as a simple form, i.e. is built out of representation
theory of $G$ over an algebracially closed field.  In general, though,
the $K$-theory of the category $V^G(F^{\prime})$ may not be a priori any
simpler than the $K$-theory of $F$.  However, when $F^{\prime}$ is a
field whose $K$-theory we already understand well, then we expect to
obtain information this way.  In addition to algebraically closed fields,
we have an understanding of the $K$-theory of finite fields from the work
of Quillen \cite{Quillen2}. So, consider  the case where
$F$ has  finite characteristic
$p$, distinct from $l$.  We may in this case let $F^{\prime}$ be the
maximal finite subfield contained in $F$, which is of the form
$\Bbb{F}_q$, where
$q = p^n$ for some $n$. In this case we have, in parallel with
Conjecture \ref{assemblyconj}, 

\begin{gConjecture} For $F$ of finite characteristic $p \neq l$, with
$\Bbb{F}_q$ the maximal finite subfield of $F$, the map 

$$K(G,\overline{\Bbb{F}_q})^{\wedge}_{{\cal H}_l} \longrightarrow
K(G,E)^{\wedge}_{{\cal H}_l}
$$
is an equivalence of $G$-$S$-algebras.  In particular, we have an
equivalence 

$$(K(G,\overline{\Bbb{F}_q})^{\wedge}_{{\cal H}_l})^G \simeq
K(G,E)^{\wedge}_{{\cal H}_l} \cong KF. 
$$
\end{gConjecture}

\noindent We now argue that this is a reasonable replacement for
Conjecture \ref{assemblyconj} in this case.  We note that the fixed point
category of $V(G,\overline{\Bbb{F}_q})$ is the category $V^G(\overline{\Bbb{F}_q})$ of linear
descent data.  We have the straightforward observation 

\begin{gProposition} For any finite group $G$ acting on $\Bbb{F}_q$ by
automorphisms, the category of linear descent data $V^G(\Bbb{F}_q)$ is
equivalent to the category of left modules over the twisted group ring 
$\Bbb{F}_q \langle G \rangle$. Similarly, for a profinite group $G$
acting on $\Bbb{F}_q$, we find that the category of linear descent data
is equivalent to the category of continuous modules over the twisted group ring  $\Bbb{F}_q \langle G \rangle$. 
\end{gProposition}

\noindent Note that for finite $G$, $\Bbb{F}_q \langle G \rangle$ is a finite dimensional
semisimple algebra over $\Bbb{F}_p$, so it is a product of matrix rings
over field extensions of $\Bbb{F}_p$.  Consequently, we should view its
$K$-theory as essentially understood, given Quillen's computations for
finite fields.  For this reason, we regard the above conjecture as a
satisfactory replacement for Conjecture \ref{assemblyconj}. 

\noindent In the characteristic 0 case, though, we have more
difficulties, since in this case we do not know the $K$-theory of the
prime field $\Bbb{Q}$.  Indeed, we would like to make conjectures about
the $K$-theory of $\Bbb{Q}$ involving the representation theory of
$G_{\Bbb{Q}}$.  We do, however, thanks to the work of Suslin
\cite{Suslin2}, understand the $K$-theory of the field $\Bbb{Q}_p$ and
$\Bbb{Z}_p$ where
$p \neq l$.  

\begin{gTheorem}\label{suslin}{\bf (Suslin; see \cite{Suslin2})}  Let $K$
denote any finite extension of $\Bbb{Q}_p$, and let ${\cal O}_K$ denote
its ring of integers. Let $(\pi)$ denote its unique maximal ideal. The
quotient {\homo}  
 ${\cal O}_K  \rightarrow {\cal O}_K /\pi {\cal O}_K$ induces an
isomorphism on
$l$-completed $K$-theory. 
\end{gTheorem} 

Now, suppose we have any field $F$ containing $\Bbb{Q}_p$ and 
the $l$-th roots of unity,   and let
$E$ denote its algebraic closure. Let
$L$ be $\bigcup _n \Bbb{Q}_p(\zeta _{l^n}) \subseteq E$, so $L \subseteq F$.  $L$ is of course
closed under the action of the absolute Galois group $G = G_F$.  By abuse
of notation, we will write $V(G,{\cal O}_L)$ for the category of twisted
$G-{\cal O}_L$-modules over
${\cal O}_L$, i.e. finitely generated ${\cal O}_L$-modules $M$ equipped
with a $G$-action so that $g(rm) = r^gg(m)$ for  all $r \in {\cal O}_L$,
$m \in M$, and $g \in G$.  The following is an easy consequence of 
Theorem \ref{suslin} above. 

\begin{gProposition} The functor $V(G, {\cal O}_L) \rightarrow V(G,
{\cal O}_L / \pi {\cal O}_L)$ induces an equivalence $K^G{\cal O}_L
\rightarrow K^G {\cal O}_L / \pi {\cal O}_L$.  
\end{gProposition}

\noindent Since ${\cal O}_L/\pi
{\cal O}_L$ is a semisimple algebra over a finite field, we will
regard it as an understood quantity.    This means that in this case, we
also have a version of the representational assembly via the diagram

$$
\begin{diagram}
\node{(K^G {\cal O}_L / \pi {\cal O}_L)^{\wedge}_{{\cal H}_l}}
 \node{ (K^G{\cal
O}_L)^{\wedge}_{{\cal H}_l}}\arrow{w,t}{\sim} \arrow{e,t}{E
\column{\otimes}{{\cal O}_L}}
\node{(K^GE)^{\wedge}_{{\cal H}_l} \cong KF}
\end{diagram}
$$

\noindent This
kind of result extends to the case where $F$ contains a Henselian local
ring whose residue class field is algebraic over a finite field.

\section{The ascent map and assembly for the case $\mu
_{l^{\infty}} \subseteq F$}\label{implicit}

\noindent The constructions of the Section \ref{secthree}  provide explicit maps
from derived completions of spectra attached to representation categories
of the absolute Galois group $G$ of a field $F$ to  $KF$. However,
they only apply to the geometric case, i.e. where $F$ contains an algebraically closed subfield.  We believe, though,  that the statement relating the $K$-theory of the field to the derived representation theory of the absolute Galois group should be true for all fields containing the $l$-th power roots of unity. The purpose of this section is to develop a criterion which when satisfied will produce an equivalence between these spectra. The criterion will involve a map which we call the
{\em ascent} map. The terminology ``ascent''
refers to the fact that the method gives a description of  $KE$ coming
from information about
$KF$, rather than describing $KF$ in terms of $KE$, as is done in the
descent.  In the interest of clarity, we will first
describe the approach as it works in the abelian case, where it is not
necessary to pass to the generality of Mackey and Green functors, and
then indicate the changes necessary when we deal with the more general
situation. Finally, we suppose that $F$ contains all the $l$-th power
roots of unity, and as usual we let $E$ denote the algebraic closure of $F$ and $G$ denote the absolute Galois group of $F$. 

Consider the forgetful functor $V^G(E) \stackrel{\phi}{\rightarrow} V(G,E)
\cong Vect(E)$, which simply forgets the $G$-action.  It induces a map 

$$ KF \cong K^G(E) \stackrel{K\phi}{\longrightarrow } K(G,E) \cong KE.
$$

\noindent We now have a commutative diagram

$$
\begin{diagram}
\node{KRep_F[G] \wedge K^G(E)} \arrow{e} \arrow{s,t}{\varepsilon \wedge
id} \node{K^G(E)} \arrow{s,b}{K\phi} \\
\node{KF \wedge K^G(E)} \arrow{e} \arrow{s,t}{id \wedge K \phi}
\node{KE} \arrow{s,b}{id} \\
\node{KF \wedge KE} \arrow{e} \node{KE}
\end{diagram}
$$

\noindent of spectra, where the horizontal maps are multiplication maps
defining algebra structures, and where the vertical maps are built out
of the augmentation map $\varepsilon$ and the forgetful functor $\phi$. 
Note that the middle horizontal map exists because $\varepsilon$ is
induced by the forgetful functor $Rep_F[G] \rightarrow Vect(F)$  which
forgets the $G$ action. Considering the four corners of this diagram, we
obtain a map of spectra

$$asc_F: KF \column{\wedge}{KRep_F[G]}K^G(E) \longrightarrow K(G,E)
\cong KE.
$$

\noindent It is readily verified that $asc_F$ is a {\homo} of commutative
$S$-algebras.  We note further that both sides of this map are $KF$-algebras, and that $asc_F$ is a homomorphism of $KF$-algebras.  As usual, let $\Bbb{H}$ denote the mod-$l$ Eilenberg-MacLane spectrum, so we have an obvious morphism $KF \rightarrow \Bbb{H}$ of $S$-algebras.  We will now explore the consequences of the hypothesis that the homomorphism

$$ id_{\Bbb{H}} \column{\wedge}{KF} KF \column{\wedge}{KRep_F[G]} K^G(E) = \Bbb{H} 
\column{\wedge}{KRep_F[G]} K^G(E) \rightarrow  
\Bbb{H} \column{\wedge}{KF} KE
$$
is a weak equivalence.

\noindent   For any homomorphism
of $S$-algebras $f  : D \rightarrow E$, we let $T^{\cdot}(D \rightarrow E
)$ denote the cosimplicial $S$-algebra given by 

$$ T^k(D \rightarrow E) = \underbrace{E \column{\wedge}{D} E
\column{\wedge}{D} \cdots \column{\wedge}{D} E}_{k+1 \mbox{ factors }}
$$
so $Tot(T^{\cdot}(D \rightarrow E ) )\cong D^{\wedge}_E$.  Now suppose that $D$ and $E$ are both $k$-algebras, where $k$ is a commutative $S$-algebra, and that we are given a homomorphism $k \rightarrow \Bbb{H}$.  We can now construct a bisimplicial spectrum $\Sigma ^{\cdot \cdot} =  \Sigma ^{\cdot \cdot} (D \rightarrow E, k \rightarrow \Bbb{H})$ by setting 

$$\Sigma ^{kl} =   \underbrace{\Bbb{H} \column{\wedge}{k} \Bbb{H}  \column{\wedge}{k} \cdots  \column{\wedge}{k} \Bbb{H}}_{k+1  \mbox{ factors }}  \column{\wedge}{k} T^l(D \rightarrow E). 
$$

\begin{gProposition}  Suppose that $k$, $D$, and $E$ are all $(-1)$-connected, and that $\pi _0 k \cong \pi _0 E \cong \Bbb{Z}$.  Then $Tot(\Sigma ^{\cdot \cdot}) \cong D^{\wedge}_{\Bbb{H} \column{\wedge}{k} E}$, where $\varepsilon$ is the composite $D \rightarrow E \rightarrow \Bbb{H} \column{\wedge}{k} E$.  If in addition, $\pi _0 D \cong \Bbb{Z}$, then  $Tot(\Sigma ^{\cdot \cdot})$ is equivalent to the $l$-adic completion of $D$. 
\end{gProposition}
\begin{Proof} The diagonal bisimplicial spectrum is canonically equivalent to the cosimplicial spectrum $ T^{\cdot}( D \rightarrow E \column{\wedge}{k} \Bbb{H})$.  By use of the K\"{u}nneth spectral sequence, it is clear that $E \column{\wedge}{k} \Bbb{H}$ is (-1) connected and that 
$\pi _0 E \column{\wedge}{k} \Bbb{H} \cong \Bbb{F}_l$.  The first  result now follows from Theorem 6.1 of \cite{completion}. The second follows immediately from Theorem 3.2, item 4 of \cite{completion}. 
\end{Proof}

\noindent Whenever we
have a commutative diagram 

$$
\begin{diagram} \node{A} \arrow{s} \arrow{e} \node{B} \arrow{s} \\
\node{A^{\prime}} \arrow{e} \node{B^{\prime}}
\end{diagram}
$$
of $k$-algebras, we have an induced map $$\Sigma^{\cdot \cdot }(A \rightarrow B)
\rightarrow \Sigma ^{\cdot \cdot }(A^{\prime} \rightarrow B^{\prime}).$$

Now suppose we have  two $k$-algebra homomorphisms $A \rightarrow
B$ and
$A
\rightarrow C$, so that each of the composites $k \rightarrow A \rightarrow B$ and $k
\rightarrow A \rightarrow C$ are weak equivalences. We then obtain a
commutative diagram of commutative $k$-algebras

$$
\begin{diagram}  \node{B} \arrow{e} \arrow{s} \node{ B
\column{\wedge}{A}C} \arrow{s} \\
\node{B \column{\wedge}{A}{C } \column{\wedge}{k}A}
\arrow{e} \node{B \column{\wedge}{A} C.}
\end{diagram}
$$
The upper horizontal arrow and the left hand vertical arrow are obvious
inclusions, the right hand arrow is the identity, and the lower
horizontal arrow is the multiplication map obtained by regarding $B
\column{\wedge}{A} C $ as a right module over the $k$-algebra $A$.  We
consequently obtain a map of bicosimplicial $k$-algebras
$$\Sigma^{\cdot \cdot }(B \rightarrow B \column{\wedge}{A} C) \longrightarrow
\Sigma^{\cdot \cdot}(B \column{\wedge}{A}{C } \column{\wedge}{k}A \rightarrow 
B \column{\wedge}{A} C). 
$$

\begin{gProposition} Under the hypotheses above, i.e. that $k \rightarrow
B$ and $k \rightarrow C$ are weak equivalences, this map of bicosimplicial
$k$-algebras is a levelwise weak equivalence, and hence induces an
equivalence  of total spectra. 
\end{gProposition}
\begin{Proof} There is a commutative diagram

$$
\begin{diagram} \node{B  \column{\wedge}{A}  T^{\cdot}(A \rightarrow C)}
\arrow{s} \arrow{e,t}{\sim} \node{T^{\cdot}(B \rightarrow  B
\column{\wedge}{A} C)} \arrow{s} \\
\node{B \column{\wedge}{A}  T^{\cdot}(C  \column{\wedge}{k} A
\rightarrow C)} \arrow{e,t}{\sim} 
\node{T^{\cdot}(B  \column{\wedge}{A}  C  \column{\wedge}{k}  A
\rightarrow B  \column{\wedge}{A}  C)}
\end{diagram}
$$
where the horizontal arrows are levelwise weak equivalences.  This follows
from the  standard isomorphism 

$$  B \column{\wedge}{A} ( \underbrace{C \column{\wedge}{A} C \column{\wedge}{A} \cdots \column{\wedge}{A} C)}_{k \mbox{ factors }} \rightarrow 
\underbrace{(B \column{\wedge}{A}C ) \column{\wedge}{B} (B \column{\wedge}{A}C ) \column{\wedge}{B}  \cdots \column{\wedge}{B} (B \column{\wedge}{A}C ) }_{k \mbox{ factors }}.
$$

It now follows  that the map $ T^{\cdot}(A \rightarrow C) \rightarrow 
T^{\cdot}(C \column{\wedge}{k} A \rightarrow C)$ is a levelwise weak
equivalence, since  the composite 
$$ A \cong k \column{\wedge}{k} A\longrightarrow C  \column{\wedge}{k} A
$$
is a weak equivalence.  Applying the functors 
$$\underbrace{\Bbb{H} \column{\wedge}{k} \Bbb{H}\column{\wedge}{k} \cdots \column{\wedge}{k} \Bbb{H} }_{k \mbox{ factors}}\column{\wedge}{k}-
$$
to this equivalence, we obtain the result for $\Sigma ^{\cdot \cdot}$.
\end{Proof}

\noindent We now examine the consequences of this result in our case. 
The spectra will be $k = KF$, $A = KRep_F[G]$, $B = K^G(E)$, and $ C =
KF$, viewed as an $A$-algebra via the augmentation which forgets the
action.  This data clearly satisfies the hypotheses above.  If 
$id_{\Bbb{H}} \column{\wedge}{KF} asc_F$ is a weak equivalence, then the arrow $B \rightarrow B
\column{\wedge}{A} C$ can be identified with the arrow $K^G(E)
\rightarrow K(G,E)$, which in turn can be identified canonically with
the arrow $KF \rightarrow KE$.  It follows that 
$$ Tot(\Sigma^{\cdot \cdot }( B \rightarrow B
\column{\wedge}{A} C) )\cong Tot(\Sigma^{\cdot \cdot }(KF \rightarrow KE)) \cong
KF^{\wedge}_l
$$
On the other hand,
$$\Sigma^{\cdot \cdot} (B  \column{\wedge}{A}  C 
\column{\wedge}{k}  A \rightarrow B  \column{\wedge}{A}  C) \cong $$
$$\Sigma^{\cdot \cdot} (KF \column{\wedge}{KRep_F[G]} K^G(E) \column{\wedge}{KF} KRep_F[G] \rightarrow KF \column{\wedge}{KRep_F[G]} K^G(E)  ) 
 $$ 
and $asc_F$ induces a map of bicosimplicial spectra
$$\Sigma^{\cdot \cdot} (KF \column{\wedge}{KRep_F[G]} K^G(E) \column{\wedge}{KF} KRep_F[G] \rightarrow KF \column{\wedge}{KRep_F[G]} K^G(E)  ) \longrightarrow
$$
$$ \Sigma ^{\cdot \cdot} (KE \column{\wedge}{KF} KRep_F[G] \rightarrow KE).
$$
By Proposition 3.2, item 2 of \cite{completion}, if $id_{\Bbb{H}} \column{\wedge}{KF} asc_F$ is a weak equivalence, this map of bicosimplicial spectra is a weak equivalence.  The natural homomorphism $KE \column{\wedge}{KF} KRep_F[G] \rightarrow KRep_E[G]$ now induces another map of bicosimplicial spectra
$$  \Sigma ^{\cdot \cdot} (KE \column{\wedge}{KF} KRep_F[G] \rightarrow KE) \longrightarrow 
\Sigma^{\cdot \cdot}(KRep_E[G] \rightarrow KE).
$$

 We now
have 
\begin{gProposition} When all $l$-th power roots of unity are in $F$, the
evident map $KRep_F[G] \column{\wedge}{KF} KE
\rightarrow KRep_E[G]$ is a weak equivalence. \end{gProposition}
\begin{Proof} It is a standard fact in representation theory that for
any field $F$ which contains the $l$-th power roots of unity, the
isomorphism classes of representations  of any finite $l$-group
are in bijective correspondence with the isomorphism classes of
representations of the same group in $\Bbb{C}$. Moreover, the
endomorphism ring of any irreducible representation is a copy of $F$.  We
therefore have isomorphisms
$$\pi _*KRep_F[G] \cong K_*F \otimes R[G]
$$
and 
$$\pi _*KRep_E[G] \cong K_*E \otimes R[G].
$$
where $R[G]$ denotes the complex representation ring.  
The $E_2$-term of the  K\"{u}nneth spectral sequence for $\pi _*
KRep_F[G] \column{\wedge}{KF} KE$  now has the form 
$$ Tor_{K_*F }(K_*F \otimes R[G], K_*E) \cong K_*Rep_E[G]$$

which gives the result.  
\end{Proof}

\noindent We now have the following conclusion.

\begin{gTheorem} \label{noneq} Suppose that $F$ contains all the $l$-th power roots of
unity.  Suppose further that $id_{\Bbb{H}} \column{\wedge}{KF} asc_F$ is a weak equivalence.  Then we
have a weak equivalence
$$ KF^{\wedge}_l \cong KRep_E[G]^{\wedge}_{KE}.
$$ 
\end{gTheorem}

\noindent We observed in Section \ref{secthree} that we did not expect to
have an equivalence 

$$ KF \cong KRep_k[G]^{\wedge}_{\Bbb{H}}
$$
as stated except in the abelian case, but that there is an analogous
``fully equivariant version'' involving equivariant spectra and Green
functors which we expect to hold in all cases.  This suggests that we
should not expect the ascent map as formulated above to hold except in
the abelian case. However, it is possible to modify the ascent map to
extend it to the fully equivariant version of the completion, and we
indicate how. 

\noindent We fix the group $G = G_F$, and we define various equivariant $S$-algebras.  We let $\underline{KV}$ denote the equivariant spectrum obtained using the category of descent data described above.  For any subgroup $H \subseteq G$, the fixed point spectrum $\underline{KV}^H$ is equivalent to the $K$-theory spectrum $K(E^H)$.  For any field $F$, we let $\underline{KRep_F}$ denote the equivariant spectrum corresponding to the trivial action of $G$ on $F$, so for any $H \subseteq G$, we have that $\underline{KRep_F}^H \cong KRep_F[H]$.  Also for any field $F$, we will write $\underline{KF}$ for the equivariant spectrum $KF \wedge S^0$, where $S^0$ denotes the $G$-equivariant sphere spectrum.  We have that $\underline{KF}^H \cong KF \wedge G/H_+$ for any subgroup of finite index $H \subseteq G$.  Finally, we let $\underline{{\cal  H}}$ denote the $G$-spectrum associated to the Green  functor 
${\cal Z}/l{\cal  Z} $ described above in section \ref{secthree}.  We note that there is an obvious homomorphism of $G$-equivariant $S$-algebras $\underline{KF} \rightarrow \underline{{\cal H}_l}$.  We can easily verify that there is a natural analogue of the ascent map in this equivariant setting, i.e. a map of $G$-equivariant $S$-algebras
$ asc_F: \underline{KF} \column{\wedge}{\underline{KRep_F}} \underline{KV} \rightarrow 
\underline{KE}$
We now have the following  equivariant generalization of \ref{noneq}, whose proof is identical to the nonequivariant version given above.

\begin{gTheorem} Let $F,E,$ and $G$ be as above.  Suppose further that the map of equivariant spectra
$$ id_{{\cal H}_l} \column{\wedge}{\underline{KF}} asc_F : {\cal H}_l \column{\wedge}{\underline{KRep_F}} \underline{KV} \rightarrow {\cal H}_l \column{\wedge}{\underline{KF}} \underline{KE}
$$
is a weak equivalence of spectra.  Then there is a canonically defined equivalence of $G$-equivariant $S$-algebras
$$  \underline{KV} ^{\wedge}_{{\cal H}_l} \rightarrow \underline{KRep_E}^{\wedge}_{{\cal H}_l}.
$$
 Note that the fixed point sets of the equivariant spectrum $\underline{KV} ^{\wedge}_{{\cal H}_l}$ are just the $l$-adic completions of the $K$-theory spectra of the corresponding subfields.  
\end{gTheorem}

{\bf Remark:} We expect that the equivariant map $ id_{{\cal H}_l} \column{\wedge}{\underline{KF}} asc_F $ will be an equivalence for all fields containing all the $l$-th power roots of unity. 

{\bf Remark:}  The above discussion could be a reflection of an interesting convergence theorem within $\Bbb{A}^1$ homotopy theory (see \cite{Morel} or \cite{MorelV} for treatments of this construction).  We suppose that $\Bbb{A}^1$ homotopy theory has been extended  to an equivariant theory for the action of  profinite groups. We further suppose that algebraic $K$-theory has been extended to a fully equivariant cohomology theory $K^G$ on this equivariant category, where by ``fully equivariant" we mean a theory which admits deloopings with respect to arbitrary  continuous linear representations of a fixed  profinite group $G$.  Let $F$ be a field, $E$ its separable closure, and suppose $G = G_F$ is its separable Galois group.  We let the group $\Gamma = G \times G$ act on $E$ by projection on the first factor.  We let $G_0 = \{ e \} \times G \subseteq \Gamma $, and let $\Delta$ denote the diagonal subgroup in $\Gamma $.   We consider the square
$$
\begin{diagram}
\node{ Spec(E) \times \Gamma } \arrow{s} \arrow{e}  \node{Spec(E) \times \Gamma / \Delta }  \arrow{s}  \\
\node{Spec(E) \times \Gamma / G_0} \arrow{e}  \node{Spec(E)}
\end{diagram} 
$$
\noindent of pro-schemes, where the sets $\Gamma$, $\Gamma / \Delta$, and $\Gamma /G_0$ are all regarded as discrete pro-shemes. If this square is a homotopy pullback square in the equivariant $\Bbb{A}^1$ category,  and if the Eilenberg-Moore spectral sequence associated to the equivariant cohomology theory $K^{\Gamma}$ converges, then this would confirm the hypotheses of the theorems above.

\section{Derived representation theory and deformation $K$-theory}\label{deformation}

\noindent So far, our understanding of the derived
completion of representation rings is limited to  knowledge of
spectral sequences for computing them, given information about $Tor$ or
$Ext$ functors of these rings.  In this (entirely speculative) section,
we want to suggest that there should be a relationship between derived
representation theory and deformations of representations.  

\noindent We consider first the representation theory of finite
$l$-groups $G$.  In this case, the derived representation ring of $G$ is
just the $l$-adic completion of $R[G]$, as can be verified using the
results of Atiyah \cite{Atiyah}. In particular, the derived representation
ring has no higher homotopy.  However, when we pass to a profinite
$l$-group, such as $\Bbb{Z}_l$, we find that the derived representation
ring does have higher homotopy, a single copy of $\Bbb{Z}_l$ in dimension
1.  This situation appears to be parallel to the following situation. 
Consider the infinite discrete group $\Gamma = \Bbb{Z}$.  We may consider
the category $Rep_{\Bbb{C}}[\Gamma]$ of finite dimensional representations
of
$\Gamma$, and its
$K$-theory spectrum $KRep_{\Bbb{C}}[\Gamma]$.  The homotopy of this
spectrum is given by 
$$ \pi _*KRep_{\Bbb{C}}[\Gamma] \cong \Bbb{Z}[S^1]\otimes K_* \Bbb{C}
$$
where $S^1$ is the circle regarded as a discrete group. This isomorphism
arises from the existence of Jordan normal form, which (suitably
interpreted) shows that every representation of $\Gamma$ admits a
filtration by subrepresentations so that the subquotients are
one-dimensional, and are therefore given by multiplication by a uniquely
defined non-zero complex number.  This construction does not take into
account the topology of $\Bbb{C}$ at all.  We might take the topology
into account as in the following definition.  

\begin{gDefinition}  Consider any discrete group $\Gamma$. For each $k$,
we consider the category $Rep^k_{\Bbb{C}}[\Gamma]$ whose objects are all
possible continuous actions of $\Gamma$ on $\Delta [k] \times V$ which
preserve the projection $\pi : \Delta [k] \times V \rightarrow \Delta
[k]$, and which are linear on each fiber of $\pi$.  It is clear that the
categories $Rep^k_{\Bbb{C}}[\Gamma]$ fit together into a simplicial
symmetric monoidal category, and we define the {\bf deformation
$K$-theory spectrum of
$\Gamma$},
$K^{def}[\Gamma]$, as the total spectrum of the simplicial spectrum $$k
\rightarrow KRep^k_{\Bbb{C}}[\Gamma].$$
\end{gDefinition}
{\bf Remark:} The terminology is justified by the observation that, for
example, an object in the category of 1-simplices
$Rep^1_{\Bbb{C}}[\Gamma]$ is exactly a path in the space of
representations of $\Gamma$ in $GL(V)$, or a deformation of the
representation at $0 \times V$. 

\noindent One can easily check that $\pi _0 K^{def}[\Bbb{Z}] \cong 
\Bbb{Z}$, with the isomorphism given by sending a represenatation to its
dimension.  This follows from the fact that any two representations of
the same dimension of 
$\Bbb{Z}$ can be connected by a deformation. In fact, one can apply the
functor
$\pi _0$ levelwise to the simplicial spectrum
$KRep^{\cdot}_{\Bbb{C}}[\Bbb{Z}]$, and attempt to compute the homotopy
groups of the simplicial abelian group 
$k \rightarrow  \pi _0 KRep^{k}_{\Bbb{C}}[\Bbb{Z}]$.  As above, it is
easy to see that $\pi _0$ of this simplicial abelian group is zero, and
it appears likely  that
$$\pi _* (\pi _0 KRep^{\cdot}_{\Bbb{C}}[\Bbb{Z}]) \cong \pi _*
\Bbb{Z}[S^1] \cong H_*(S^1)
$$ 
where $S^1$ is the circle regarded as a topological group. Note that $S^1$
is the character  group of the original discrete group $\Bbb{Z}$. 

\noindent Now consider the profinite group $\Bbb{Z}_l$. We may define the
deformation $K$-theory of this group as the direct limit of its finite
quotients, and in view of the rigidity of complex representations of
finite groups we have that
$$KRep_{\Bbb{C}}[\Bbb{Z}_l] \cong K^{def}[\Bbb{Z}_l.]
$$
 We have an inclusion
$\Bbb{Z}
\hookrightarrow
\Bbb{Z}_l$, and  therefore a map of spectra $K^{def}[\Bbb{Z}_l]
\rightarrow  K^{def}[\Bbb{Z}]$. We now obtain a composite map of derived
completions 
$$KRep_{\Bbb{C}}[\Bbb{Z}_l]^{\wedge}_{{\Bbb{ H}}} \simeq
K^{def}[\Bbb{Z}_l]^{\wedge}_{{\Bbb{ H}}}
\longrightarrow K^{def}[\Bbb{Z}]^{\wedge}_{\Bbb{ H}}
$$
 where $\Bbb{H}$ denotes the mod-$l$ Eilenberg-MacLane spectrum, regarded as a module over the deformation $K$-theory via the map which takes any representation to its dimension mod $l$.  
Since the homotopy groups on the two sides appear isomorphic, it appears
likely that this composite is a homotopy equivalence.  We now observe
that for any discrete group $\Gamma$, we may consider its pro-$l$
completion $\Gamma ^{\wedge}_l$, and obtain a map $KRep_{\Bbb{C}}[\Gamma
^{\wedge}_l] \rightarrow K^{def}[\Gamma]$.  This map induces a map $j$ of
derived completions, and we ask the following question. 

\noindent {\bf Question:} For what discrete groups does the map $j$
induce a weak equivalence of spectra after derived completion at the
augmentation map $\varepsilon : K^{def}[\Gamma] \rightarrow  K\Bbb{C}$?
Note that the completion may have to be taken in the category of
$\Gamma^{\wedge}_l$-equivariant spectra. 

\noindent We can further ask if the effect of derived completion can be
computed directly on the pro-$l$ group, rather than by permitting
deformations on a discrete subgroup.  This cannot be achieved over
$\Bbb{C}$ since representations of finite groups are rigid, but perhaps
over some other algebraically closed field. 

\noindent {\bf Question:} Can one construct a deformation $K$-theory of
representations of a pro-$l$ group $G$, using an algebraic deformation
theory like the one discussed in \cite{Mazur}, so that it coincides with
the derived completion of $KRep_{\Bbb{C}}[G]$? 

\noindent We conclude by pointing out an analogy between our theory of
the representational assembly and other well known constructions. 
Consider a $K(\Gamma, 1)$ manifold $X$.  The category $VB(X)$ of complex
vector bundles over $X$ is a symmetric monoidal category, and we can
construct its $K$-theory. On the other hand, we have a functor from the
category of complex representations of $\Gamma$ to $VB(X)$, given by 
$$ \rho \rightarrow (\tilde{X} \column{\times}{\Gamma} \rho \rightarrow X)
$$ 
where $\tilde{X}$ denotes the universal covering space of $X$. This
functor produces a map of $K$-theory spectra.  It is also easy to define
a ``deformation version'' of $KVB(X)$, and one obtains a map 
$$K^{def}[\Gamma] \rightarrow K^{def}VB(X).
$$
This map should be viewed as the analogue in this setting for the
representational assembly we discussed in the case of $K$-theory of
fields, using the point of view that fields are analogues of $K(G,1)$
manifolds in the algebraic geometric context.

\end{document}